\documentclass[final]{siamart}
\usepackage{mwe}
\usepackage{subfig}

\usepackage{lipsum}
\usepackage{amsfonts}
\usepackage{graphicx}
\usepackage{epstopdf}
\usepackage{algorithmic}
\ifpdf
  \DeclareGraphicsExtensions{.eps,.pdf,.png,.jpg}
\else
  \DeclareGraphicsExtensions{.eps}
\fi

\newcommand{\TheTitle}{Structure Preserving Model Reduction of Parametric Hamiltonian Systems} 
\newcommand{\TheAuthors}{B. Maboudi Afkham, and J. S. Hesthaven}

\headers{Symplectic Model Reduction of Hamiltonian Systems}{\TheAuthors}

\title{{\TheTitle}}

\author{
	Babak Maboudi Afkham\thanks{Department of Mathematics, Chair of Computational Mathematics and Simulation Science (MCSS), \'Ecole Polytechnique F\'ed\'erale de Lausanne, Switzerland (\email{babak.maboudi@epfl.ch})} 
	\and 
	Jan S. Hesthaven\thanks{Department of Mathematics, Chair of Computational Mathematics and Simulation Schience (MCSS), \'Ecole Polytechnique F\'ed\'erale de Lausanne, Switzerland (\email{jan.hesthaven@epfl.ch})}}

\usepackage{amsopn}

\ifpdf
\hypersetup{
  pdftitle={\TheTitle},
  pdfauthor={\TheAuthors}
}
\fi

\newcommand{\edit}{\color{black}}

\begin{document}

\maketitle

\begin{abstract}
While reduced-order models (ROMs) {\edit have been} popular for {\edit efficiently} solving large systems of differential equations, the stability of reduced models over long-time integration {\edit is of present challenges}. We present a greedy approach for ROM generation of parametric Hamiltonian systems {\edit that} captures the symplectic structure of Hamiltonian systems to ensure stability of the reduced model. Through the greedy selection of basis vectors, {\edit two new vectors are added at each iteration to the linear vector space} to increase the accuracy of the {\edit reduced} basis. We {\edit use} the error in the Hamiltonian due to model reduction as an error indicator to search the parameter space and identify the next best basis vectors. Under natural assumptions on the set of all solutions of the Hamiltonian system under variation of the parameters, we show that the greedy algorithm converges with exponential rate. Moreover, we demonstrate that combining the greedy basis with the discrete empirical interpolation method also preserves the symplectic structure. This enables the reduction of the computational cost for nonlinear Hamiltonian systems. The efficiency, accuracy, and stability of this model reduction technique is illustrated through simulations of the parametric wave equation and the parametric Schr\"odinger equation.
\end{abstract}



\section{Introduction}
{\edit Parameterized partial differential equations often arise as a model in many problems in engineering and the applied sciences}. While the need for more accuracy has led to the development of exceedingly complex models, the {\edit limitations in computational cost and storage} often make direct approaches {\edit impractical}. Hence, we must seek alternative methods that allow us to approximate the desired output under variation of the input parameters while keeping the computational costs to a minimum.

Reduced basis methods have emerged as a powerful approach for the reduction of the intrinsic complexity of such models \cite{Ito:1998up,Ito:1998ch,Ito:2001ev,Peterson:1989ki}. These methods contain two stages: {\edit the offline stage and the online stage}. In the offline stage, one explores the parameter space to construct a low-dimensional basis that accurately represents the parametrized solution to the partial differential equation. In this stage, the evaluation of the solution of the original model for multiple parameter values is required. The online stage comprises a Galerkin projection onto the span of the reduced basis, which allows exploration of the parameter space at a significantly reduced complexity \cite{Antoulas:2005:ALD:1088857,Anonymous:2016wl}.

Convectional reduced basis techniques, such as proper orthogonal decomposition (POD) \cite{Kunisch:2002er,Atwell:2001ja,Ravindran:2002hn}, require the exploration of the entire parameter space. This leads to a very expensive and often impractical offline {\edit stage} when dealing with multi-dimensional parameter domains. On the other hand, sampling techniques, usually of a greedy nature, search through the parameter space selectively, guided by an error estimate to certify the accuracy of the basis. This approach, accompanied with an efficient sampling procedure, {\edit balances the cost} of computation with the overall accuracy of the reduced-basis \cite{Cuong:2005gd,Rozza:2005ie,Anonymous:2016wl}.

{\edit Besides} computational complexity, another aspect of reduced order modeling is the preservation of structure and, in particular, {\edit the} stability of the original model. In general, reduced order models do not guarantee that such properties are preserved \cite{Anonymous:pMn0O0Q4}. 

In the context of Hamiltonian and Lagrangian systems, recent work suggests modifications of POD to preserve {\edit some} geometric structures. Lall et al. \cite{Lall:2003iy} and Carlberg et al. \cite{Carlberg:2014ky} suggests that the reduced-order system should be identified by a Lagrangian function on a low-dimensional configuration space. In this way, the geometric structure of the original system is inherited by the reduced system. {\edit Model reduction for port-Hamiltonian systems can be found in the works of Beattie et al. \cite{Chaturantabut:2016he}, Polyuga et al. \cite{Polyuga:2010gj} and references therein. These works construct a reduced port-Hamiltonian system using Krylov or POD methods that inherit the passivity and stability of the original system.} For Hamiltonian systems, Peng et al. \cite{Peng:2014di}, using a symplectic transformation, constructs a reduced Hamiltonian, as an approximation to the Hamiltonian of the original system. As a result, the reduced system preserves the symplectic structure. Although these methods preserve {\edit the} geometric structure, they use a POD-like approach for constructing the reduced basis and are not well {\edit suited} for problems with a high-dimensional parameter domain.

In this paper, we present a greedy approach for the construction of a reduced system that preserves the geometric structure of Hamiltonian systems. This technique results in a reduced Hamiltonian system that mimics the symplectic properties of the original system and preserves the Hamiltonian structure and its stability over the course of time. On the other hand, since time integration of the original system is only required once per iteration, the proposed method saves substantial computational cost during the offline stage when compared to alternative POD-like approaches. {\edit It is well known that structured matrices, e.g. symplectic matrices, generally are not well-conditioned \cite{Karow:2006cf}. The greedy update of the symplectic basis presented here, yields a orthosymplectic basis and, therefore, a norm bounded basis.} Moreover, we demonstrate that assumptions, natural for the set of all solutions of the original Hamiltonian system under variation of parameters, lead to exponentially fast convergence of the greedy algorithm. For nonlinear Hamiltonian systems, we show how the basis can be combined with the discrete empirical interpolation method (DEIM) {\edit \cite{Chaturantabut:2010cz,Barrault:2004kz}} to enable a fast evaluation of nonlinear terms while maintaining the symplectic structure.

This paper is organized as follows. Section \ref{chap:MoOr:1} presents a brief overview of model order reduction, POD and DEIM. In Section \ref{chap:Hasy:1} we cover the required topics from symplectic geometry and Hamiltonian systems. Section \ref{chap:SyMo:1} discusses the greedy generation of a symplectic reduced basis as well as other SVD-based symplectic model reduction techniques. Accuracy, stability, and efficiency of the greedy method compared to other SVD-based methods are discussed in Section \ref{chap:NuRe:1}. {\edit Finally we offer some conclusive remarks in Section \ref{chap:Con:1}}.

\section{Model Order Reduction} \label{chap:MoOr:1}
Consider a parameterized, finite dimensional dynamical system described by a set of first order ordinary differential equations
\begin{equation} \label{eq:MoOr:1}
\left\{
\begin{split}
& \frac{d}{dt} \mathbf{x}(t,\omega) = \mathbf f (t,\mathbf x,\omega), \\
& \mathbf x(0,\omega) = \mathbf x_0(\omega).
\end{split}
\right.
\end{equation}
Here $\mathbf x \in \mathbb R^n$ is the state vector, $\omega \in \Gamma$ is a vector containing all the parameters of the system {\edit belonging} to a compact set $\Gamma$ ($\subset \mathbb R^d$) and $\mathbf f : \mathbb R \times \mathbb R^n \times \Gamma \to \mathbb R^n$ is a general {\edit vector valued} function of the state variables and parameters.

We define the solution manifold as the set of all solutions to (\ref{eq:MoOr:1}) under variation of {\edit the parameters in $\Gamma$}
\begin{equation} \label{eq:MoOr:2}
	\edit \mathcal M = \{ \mathbf x(t,\omega) | \omega \in \Gamma,\ t\geq 0 \} \subset \mathbb R^n.
\end{equation}
Note that the exact solution and solution manifold is {\edit often not} available; we assume that we have a numerical integrator that can approximate the solution to (\ref{eq:MoOr:1}) for any realization of $\omega$ {\edit with a given accuracy}. By abuse of notation, we refer to $\mathbf x$ and $\mathcal M$ as the exact solution and the exact solution manifold, respectively, rather than the discrete solution and discrete solution manifold. 

Model order reduction is based on the assumption that $\mathcal M$ is of low dimension \cite{Anonymous:2016wl,Antoulas:2005:ALD:1088857} and that the span of appropriately chosen basis vectors $\{v_i\}_{i=1}^k$ covers most of the solution manifold {\edit to within} a small error. The set $\{v_i\}_{i=1}^k$ is denoted as the \emph{reduced basis} and its span as the \emph{reduced space}. Assuming that a $k$-dimensional $(k\ll n)$ reduced basis is given, the approximated solution can be represented as
\begin{equation} \label{eq:MoOr:3}
	\mathbf x \approx V \mathbf y,
\end{equation}
where $V$ is a matrix containing the reduced basis vectors as its columns and $\mathbf y$ {\edit contains} the coordinates of the approximation in this basis. By substituting (\ref{eq:MoOr:3}) into (\ref{eq:MoOr:1}) we obtain the overdetermined system
\begin{equation} \label{eq:MoOr:4}
	V \frac{d}{dt} \mathbf y = \mathbf f (t , V \mathbf y , \omega) + \mathbf r(t,\omega).
\end{equation}
Here we added the residual $\mathbf r$ to emphasize that (\ref{eq:MoOr:4}) is an approximation of (\ref{eq:MoOr:1}). Taking the Petrov-Galerkin projection \cite{Antoulas:2005:ALD:1088857} we construct a basis $W$ of size $n-k$ that is orthogonal to the residual $\mathbf r$ and {\edit requires that} $W^T V$ is invertible. This yields
\begin{equation} \label{eq:MoOr:5}
	\frac{d}{dt} \mathbf y = (W^TV)^{-1} \mathbf f(t,V\mathbf y,\omega).
\end{equation}
Equation (\ref{eq:MoOr:5}) consists of $k$ equations and is called the reduced system. Solving the reduced system instead of the original system can reduce the computational costs {\edit provided} $k$ is significantly smaller than $n$. For nonlinear systems, the evaluation of $\mathbf f$ may still have computational complexity that depends on $n$. We return to this question in detail in Section \ref{chap:MoOr.DEIM:1}.

\subsection{Proper Orthogonal Decomposition} \label{chap:MoOr.PrOr:1}
Let $\mathbf x (t_i,\omega_j)$ with $i=1,\dots,m$ and $j=1,\dots,n$ be a finite number of samples, referred to as \emph{snapshots}, from the solution manifold (\ref{eq:MoOr:2}). If we {\edit assume} that a reduced basis $V$ is provided, the projection operator from $\mathbb R^n$ onto the reduced space can be constructed as $VV^T$. The proper orthogonal decomposition (POD) requires the total error of projecting all the snapshots onto the reduced space to be {\edit minimized}. The POD basis of size $k$ is thus the solution to the optimization problem
\begin{equation} \label{eq:MoOr:6}
\begin{aligned}
& \underset{V\in \mathbb R^{n\times k}}{\text{minimize}}
& & \| S - VV^TS\|_F \\
& \text{subject to}
& & V^TV = I_k
\end{aligned}
\end{equation}
Here $S$ is the snapshot matrix, containing snapshots $\mathbf x(t_i,\omega_j)$ in its columns, $\|\cdot \|_F$ is the Frobenius norm and $I_k$ is the identity matrix of size $k$. According to {\edit Schmidt-Mirsky-Eckart-Young theorem \cite{Markovsky:2011:LRA:2103589}}, the solution to (\ref{eq:MoOr:6}) is equivalent to the truncated singular value decomposition (SVD) of the snapshot matrix $S$ given by
\begin{equation} \label{eq:MoOr:7}
	V = \sigma_1 u_1 v^T_1 + \cdots + \sigma_k u_k v^T_k.
\end{equation}
Here $\sigma_i, u_i$ and $v_i$ are the singular values, the left singular vectors, and the right singular vectors of $S$, respectively {\edit \cite{Markovsky:2011:LRA:2103589} }.

\subsection{Discrete Empirical Interpolation Method (DEIM)} \label{chap:MoOr.DEIM:1} \nocite{Chaturantabut:2010cz}
In this section we {\edit discuss the efficiency} of evaluating nonlinearities {\edit in the context of} projection based reduced models. Suppose that the right hand side in (\ref{eq:MoOr:1}) is of the form $\mathbf f(t,\mathbf x , \omega) = L\mathbf x + \mathbf g(t,\mathbf x ,\omega)$, where $L\in \mathbb R^{n\times n}$ reflects the linear part, and $\mathbf g$ is a nonlinear function. Now {\edit assume that} a $k$-dimensional reduced basis $V$ is provided. The reduced system takes the form
\begin{equation} \label{eq:MoOr:8}
	\frac{d}{dt} \mathbf y = \underbrace{(WV)^{-1} L V}_{\tilde L} \mathbf{y} + \underbrace{(WV)^{-1} \mathbf g(t,V\mathbf y,\omega)}_{\tilde N (\mathbf y)}.
\end{equation}
Here, $\tilde L$ is a $k\times k$ matrix which can be computed before time integration of the reduced system. However, {\edit the} evaluation of $\tilde N (\mathbf y)$ has a complexity that depends on $n$, the size of the original system. Suppose that the evaluation of $\mathbf g$ with $n$ components has the complexity $\alpha(n)$, for some function $\alpha$. Then the complexity of evaluating $\tilde N(\mathbf y)$ is $\mathcal{O}(\alpha(n) + 4nk)$ which consists of 2 matrix-vector operations and the {\edit evaluation of the nonlinear function, i.e. the evaluation of the nonlinear terms can be as expensive as solving the original system.}

{\edit To overcome this bottleneck we take an approach similar to that of Section \ref{chap:MoOr.PrOr:1} \cite{Chaturantabut:2010cz,Barrault:2004kz}}. Assume that the manifold $\mathcal M_{\mathbf g} = \{ \mathbf g(t,\mathbf x , \omega)| t\in \mathbb R, \mathbf x \in \mathbb R , \omega \in \Gamma\}$ is of a low dimension and that $\mathbf g$ can be approximated by a linear subspace of dimension $m\ll n$, spanned by the basis $\{ u_1 , \dots , u_m \}$, i.e.
\begin{equation} \label{eq:MoOr:10}
	\mathbf g(t,\mathbf x,\omega) \approx U \mathbf c(t,\mathbf x,\omega).
\end{equation}
Here $U$ contains basis vectors $u_i$ and $\mathbf c$ is the vector of coefficients. Now suppose $p_1,\dots,p_m$ are $m$ indices from $\{1,\dots,n\}$ and define an $n\times m$ matrix
\begin{equation} \label{eq:MoOr:11}
	P = [e_{p_1},\dots,e_{p_m}],
\end{equation}
where $e_{p_i}$ is the $p_i$-th column of the identity matrix $I_n$. Multiplying $P$ with $\mathbf g$ selects components $p_1,\dots,p_m$ of $\mathbf g$. If we assume that $P^TU$ is non-singular, the coefficient vector $\mathbf c$ can be uniquely determined from
\begin{equation} \label{eq:MoOr:12}
	P^T \mathbf g = (P^TU)\mathbf c.
\end{equation}
Finally the approximation of $\mathbf g$ is determined by
\begin{equation} \label{eq:MoOr:13}
	\mathbf g(t,\mathbf x,\omega) \approx U \mathbf c(t,\mathbf x,\omega) = U (P^TU)^{-1} P^T \mathbf g(t,\mathbf x,\omega),
\end{equation}
which is referred to as the \emph{Discrete Empirical Interpolation} (DEIM) approximation \cite{Chaturantabut:2010cz}. Applying DEIM to the reduced system (\ref{eq:MoOr:5}) yields
\begin{equation} \label{eq:MoOr:14}
	\frac{d}{dt} \mathbf y = \tilde L \mathbf y + (WV)^{-1} U(P^TU)^{-1}P^T \mathbf g(t,V\mathbf y , \omega).
\end{equation}
Note that the matrix $(WV)^{-1} U(P^TU)^{-1}$ can be computed offline and since $\mathbf g$ is evaluated only at $m$ of its components, the evaluation of the nonlinear term in (\ref{eq:MoOr:14}) does not depend on $n$.

{\edit To obtain the projection basis $U$, the POD can be applied to the ensemble of samples of the nonlinear term $\mathbf g(t_i,\mathbf x, \omega_j)$ with $i=1,\dots,m$ and $j=1,\dots,n$}. There is no additional cost {\edit associated with} computing the nonlinear snapshots, since they are generated when computing the trajectory snapshot matrix $S$. The interpolating indices $p_1,\dots,p_m$ can be constructed as follows. Given the projection basis $U = \{u_1,\dots,u_m\}$, the first interpolation index $p_1$ is chosen according to the component of $u_1$ with the largest magnitude. The rest of the interpolation indices, $p_2,\dots,p_m$ correspond to the component of the largest magnitude of the residual vector $\mathbf r = u_l - U \mathbf c$. It is shown in \cite{Chaturantabut:2010cz} that if the residual vector is a nonzero vector in each iteration then $P^TU$ is non-singular and (\ref{eq:MoOr:13}) is well defined. 

\begin{algorithm} 
\caption{Discrete Empirical Interpolation Method} \label{alg:MoOr:1}
{\bf Input:}  Basis vectors $\{u_1,\dots , u_m\}\subset \mathbb R^n$
\begin{enumerate}
\item pick $p_1$ to be the index of the largest component of $u_1$.
\item $U \leftarrow [u_1]$
\item $P \leftarrow [p_1]$
\item \textbf{for} $i\leftarrow 2$ \textbf{to} $m$
\item \hspace{0.5cm} solve $(P^TU)\mathbf c = P^T u_i$ for $\mathbf c$
\item \hspace{0.5cm} $\mathbf r \leftarrow u_i - U\mathbf c$
\item \hspace{0.5cm} pick $p_i$ to be the index of the largest component of $\mathbf r$
\item \hspace{0.5cm} $U \leftarrow [u_1,\dots,u_i]$
\item \hspace{0.5cm} $P \leftarrow [p_1,\dots,p_i]$
\item \textbf{end for}
\end{enumerate}
\vspace{0.5cm}
{\bf Output:} Interpolating indices $\{p_1,\dots,p_m\}$
\end{algorithm}

{\edit The numerical solution of (\ref{eq:MoOr:8}) may involve the computation of the Jacobian of the nonlinear function $\mathbf g(t,\mathbf x, \omega)$ with respect to the reduced state variable $\mathbf y$}
\begin{equation} \label{eq:MoOr:9}
	\edit \mathbf J_{\mathbf y}(\mathbf g) = (WV)^{-1} \mathbf J_{\mathbf x}(\mathbf g) V,
\end{equation}
{\edit where $\mathbf J_\alpha(\mathbf g)$ is the Jacobian matrix of $\mathbf g$ with respect to the variable $\alpha$.} The complexity of (\ref{eq:MoOr:9}) is $\mathcal{O}(\alpha(n) +2n^2k+2nk^2+2nk)$, comprising several matrix-vector multiplications and an evaluation of the Jacobian which depends on the size of the original system. Approximating the Jacobian in (\ref{eq:MoOr:9}) is usually both problem and discretization dependent. Often the nonlinear function $\mathbf g$ is evaluated component-wise i.e.
\begin{equation} \label{eq:MoOr:15}
	\mathbf g(\mathbf x) =
	\begin{pmatrix}
		g_1(x_1,\dots,x_n) \\
		g_2(x_1,\dots,x_n) \\
		\vdots \\
		g_n(x_1,\dots,x_n)
	\end{pmatrix}
	=
	\begin{pmatrix}
		g_1(x_1) \\
		g_2(x_2) \\
		\vdots \\
		g_n(x_n)
	\end{pmatrix}.
\end{equation}
In such cases the interpolating index matrix $P$ and the nonlinear function $\mathbf g$ commute, i.e.,
\begin{equation} \label{eq:MoOr:16}
	\tilde N(\mathbf y) \approx (WV)^{-1} U(P^TU)^{-1}P^T \mathbf g(V\mathbf y) = (WV)^{-1} U(P^TU)^{-1}\mathbf g(P^TV\mathbf y)
\end{equation}
If we now take the Jacobian of the approximate function we recover
\begin{equation}
	\edit \mathbf J_{\mathbf y}(\mathbf g) = \underbrace{ (WV)^{-1} U(P^TU)^{-1} }_{k\times m} \underbrace{ \mathbf J_{\mathbf x}(\mathbf g(P^T V \mathbf y) ) }_{m\times m} \underbrace{P^T V}_{m\times k}.
\end{equation}
The matrix $(WV)^{-1} U(P^TU)^{-1}$ can be computed offline and the Jacobian is evaluated only for $m\times m$ components. Hence the overall complexity of computing the Jacobian is now independent of $n$. {\edit We refer the reader to \cite{Barrault:2004kz,Chaturantabut:2010cz} for more detail.}

\section{Hamiltonian Systems and Symplectic Geometry} \label{chap:Hasy:1}
Let $\edit \mathcal M$ be a manifold and $\edit \Omega:\mathcal M \times \mathcal M \to \mathbb R$ be a closed, nondegenerate {\edit and skew-symmetric} 2-form on $\edit \mathcal M$. The pair $\edit (\mathcal M,\Omega)$ is called a \emph{symplectic manifold} {\edit \cite{Marsden:1999ck}}.

{\edit Let $(\mathcal M,\Omega)$ be a symplectic manifold and suppose that $H:\mathcal M \to \mathbb R$ is a smooth scalar function. The differential of $H$, denoted by $\mathbf dH$, defines a 1-form on $\mathcal M$. The nondegeneracy of $\Omega$ implies that there is a unique vector field $X_H$, \emph{Hamiltonian vector field }\cite{da2003introduction,Marsden:1999ck}, on $\mathcal M$ such that}
\begin{equation} \label{eq:Hasy:1}
	\edit i_{X_H} \Omega = \mathbf dH, 
\end{equation}
{\edit where $i_{X_H} \Omega$ is the interior product of $X_H$ with $\Omega$, i.e., that requiring}
\begin{equation}
	\edit \Omega(X_H,Y) = \mathbf dH(Y),
\end{equation}
{\edit for any vector field $Y$ on $\mathcal M$.} Note that when $\edit \mathcal M$ belongs to a Euclidean space then $\mathbf d H = \nabla_z H$. The equations of evolution are then defined by
\begin{equation} \label{eq:Hasy:2}
	\dot z = X_H(z)
\end{equation}
and known as \emph{Hamilton's equation} \cite{Marsden:1999ck}. A fundamental feature of Hamiltonian systems is the conservation of the Hamiltonian along integral curves on $\edit \mathcal M$. To emphasize the importance of this property we recall {\edit \cite{Marsden:1999ck}}

\begin{theorem} \label{theorem:Hasy:1}
Suppose that $X_H$ is a Hamiltonian vector field with the flow $\phi_t$ on a symplectic manifold $\mathcal M$. Then $H\circ \phi_t = H$.
\end{theorem}

\begin{proof}
$H$ is constant along integral curves since
\begin{equation} \label{eq:Hasy:3}
\begin{aligned}
	\frac{d}{dt}(H\circ \phi_t)(z) &= \mathbf d H(\phi_t(z)) \cdot( \frac{d}{dt} \phi_t(z) ) \\
	&= \mathbf d H (\phi_t(z))\cdot X_H(\phi_t(z)) \\
	&= \Omega_z( X_H(\phi_t(z)), X_H(\phi_t(z)) ) = 0,
\end{aligned}
\end{equation}
{\edit by using} the chain rule and bilinearity of $\Omega$ in the argument.
\end{proof}

For the case where the symplectic manifold is also a linear vector space, the pair $({\edit \mathcal M},\Omega)$ is also referred to as a \emph{symplectic vector space}. We {\edit need} the following theorems regarding symplectic vector spaces and refer the reader to \cite{de2006symplectic,Marsden:1999ck,de2006symplectic,Silva01lectureson} for detailed proofs.

{\edit
\begin{theorem} \label{theorem:Hasy:1.1} \cite{Marsden:1999ck}
If $(V,\Omega)$ is a symplectic vector space then $\Omega$ is a constant form, that is $\Omega_z$ is independent of $z\in V$. 
\end{theorem}
\begin{theorem} \label{theorem:Hasy:1.2} \cite{Marsden:1999ck}
If $(V,\Omega)$ is a finite-dimensional symplectic manifold then $V$ is even dimensional.
\end{theorem}
\begin{theorem} \label{theorem:Hasy:1.3} \cite{de2006symplectic}
(The Symplectic Gram-Schmidt) If $(V,\Omega)$ is a $2n$-dimensional symplectic vector space, then there is a basis $e_1,\dots e_n,f_1, \dots , f_n$ of $V$ such that
\begin{equation} \label{eq:Hasy:4}
\begin{aligned}
	& \Omega(e_i,e_j) = 0 = \Omega(f_i,f_j), \quad & i\neq j,\\
	& \Omega(e_i,f_j) = \delta_{ij}, & i\leq i,j \leq n.
\end{aligned}
\end{equation}
where $\delta$ is the Kronecker's delta function. Moreover if $V = \mathbb{R}^{2n}$ we can choose basis vectors $\{e_i,f_i\}_{i=1}^n$ such that
\begin{equation} \label{eq:Hasy:5}
	\Omega(v_1,v_2) = v_1^T \mathbb J_{2n} v_2, \qquad v_1,v_2\in \mathbb R^n,
\end{equation}
with $\mathbb J_{2n}$ being the symplectic matrix, defined as
\begin{equation} \label{eq:Hasy:6}
	\mathbb{J}_{2n} = 
	\begin{pmatrix}
		0_n & I_n \\
		-I_n & 0_n
	\end{pmatrix}.
\end{equation}
Here $I_n$ and $0_n$ is the identity matrix and the zero square matrix of size $n$, respectively.
\end{theorem}
\begin{theorem} \label{theorem:Hasy:1.4} \cite{Marsden:1999ck}
The classical inner product $\langle \cdot,\cdot \rangle:\mathbb R^{2n}\times \mathbb R^{2n}\to \mathbb R$ can be written in terms of the 2-form as
\begin{equation}
	\langle v,u \rangle = \Omega(\mathbb J_{2n}v,u),\quad \forall u,v \in \mathbb R^{2n}.
\end{equation}
\end{theorem}
}

{\edit 
\begin{definition}  \cite{de2006symplectic}
Suppose $(V,\Omega)$ is a finite dimensional symplectic vector space and $E\subset V$ is a subspace. Then the symplectic complement of $E$ inside $V$ is defined as
\[
	E^{\perp} := \{ v\in V |\ \Omega(v,e) = 0,\ \forall e\in E \}
\]
\end{definition}
Note that $E \cap E^{\perp}$ is not empty in general. 
\begin{definition} \cite{de2006symplectic}
Suppose $(V,\Omega)$ is a finite dimensional symplectic vector space. A subspace $E\subset V$ is called a Lagrangian subspace inside $V$ if $E = E^\perp$.
\end{definition}
\begin{theorem} \label{theorem:Hasy:1.5} \cite{abraham1978foundations}
Suppose $(V,\Omega)$ is a finite dimensional symplectic vector space. If $E\subset V$ is a Lagrangian subspace then $dim(E)=\frac 1 2dim(V)$. Here $dim$ denotes the dimension of the subspace.
\end{theorem}
\begin{definition}
A basis of $(V,\Omega)$ is called orthosymplectic if it is both a symplectic basis and an orthogonal basis with respect to the classical scalar product.
\end{definition}
\begin{theorem} \label{theorem:Hasy:1.6}  \cite{da2003introduction}
Suppose $(V,\Omega)$ is a $2n$ dimensional symplectic vector space and $E\subset V$ is a Lagrangian subspace. Then there is an orthosymplectic basis for $V$.
\end{theorem}
\begin{proof}
Starting from a Lagrangain subspace in $E \subset V$ an orthosymplectic basis can be easily constructed. By Theorem \ref{theorem:Hasy:1.5} $E$ is $n$ dimensional. Suppose that $\{ e'_1,\dots, e'_n \}$ is a basis for $E$, using the classical Gram-Schmidt orthogonalization process we can construct an orthonormal basis $\{ e_1,\dots,e_n \}$. Define a new set of vectors $f_1 = \mathbb J_{2n}^Te_1$, $f_2 =\mathbb J_{2n}^T e_2$, $\dots$, $f_n= \mathbb J_{2n}^Te_n$. We have
\begin{equation}
	\langle f_i, f_j \rangle = e_i^T \mathbb J_{2n} {\mathbb J_{2n}}^T e_i = \delta_{ij}, \quad \langle f_i, e_j \rangle = e_i^T \mathbb J_{2n} e_j = 0, \quad i,j=1,\dots,n,
\end{equation}
where we used the fact that $\mathbb J_{2n} {\mathbb J_{2n}}^T = I_{2n}$ in the first identity and the second identity is due to the fact that the basis $\{ e_1,\dots,e_n \}$ forms a Lagrangian subspace. This shows that the set $\{ e_1,\dots,e_n \}\cup \{ f_1,\dots,f_n \}$ forms an orthonormal basis. Also, it can be easily verified that this is a symplectic basis. Thus $\{ e_1,\dots,e_n \}\cup \{ f_1,\dots,f_n \}$ is an orthosymplectic basis.
\end{proof}
\begin{theorem} \label{theorem:Hasy:1.7} \cite{Marsden:1999ck}
On a finite-dimensional symplectic vector space the relationship (\ref{eq:Hasy:1}) becomes 
\begin{equation} \label{eq:Hasy:7}
\left\{
\begin{aligned}
	&\dot {\mathbf z} = \mathbb{J}_{2n} \nabla_{\mathbf z} H(\mathbf z), \\
	& \mathbf z(0) = \mathbf z_0.
\end{aligned}
\right.
\end{equation}
or, by introducing the canonical coordinates $\mathbf z = (\mathbf q^T, \mathbf p^T)^T$,
\begin{equation} \label{eq:Hasy:8}
\left\{
\begin{aligned}
	&\dot {\mathbf q} = \nabla_{\mathbf p} H(\mathbf q,\mathbf p),\\
	&\dot {\mathbf p} = - \nabla_{\mathbf q} H(\mathbf q,\mathbf p).
\end{aligned}
\right.
\end{equation}
\end{theorem}
}

{\edit Let} us now introduce \emph{symplectic transformations}, i.e., mappings between symplectic manifolds which preserve the 2-form $\Omega$. The accurate numerical treatment of Hamiltonian systems often requires preservation of the symmetry expressed in Theorem \ref{theorem:Hasy:1}. Symplectic transformations can be used to construct such symmetry preserving numerical methods. 

{\edit
\begin{definition}
Let $(V,\Omega)$ and $(W,\Pi)$ be two linear symplectic vector spaces of dimensions $2n$ and $2k$, respectively. A linear mapping $\phi:V \to W$ is called \emph{symplectic} or \emph{canonical} if
\begin{equation} \label{eq:Hasy:9}
	\Omega = \phi^* \Pi
\end{equation}
where $\phi^* \Pi$ is the pullback of $\Pi$ by $\phi$, i.e. for all $\mathbf{z}_1, \mathbf{z}_2\in V$
\begin{equation}
	\Omega(\mathbf{z}_1,\mathbf{z}_2) = \Pi(\phi(\mathbf{z}_1),\phi(\mathbf{z}_2)).
\end{equation}
\end{definition}

Note that if we represent the transformation $\phi$ as a matrix $A\in \mathbb R^{2n\times 2k}$ condition (\ref{eq:Hasy:9}) is equivalent to \cite{Marsden:1999ck}}

\begin{equation} \label{eq:Hasy:11}
	A^T \mathbb{J}_{2n}A = \mathbb{J}_{2k}.
\end{equation}
A matrix of size $2n\times 2k$ satisfying (\ref{eq:Hasy:11}) is called a \emph{symplectic matrix}.

\begin{definition}
	The \emph{symplectic inverse} of a matrix $A\in \mathbb{R}^{2n\times 2k}$ is denoted by $A^+$ and defined by {\edit \cite{Peng:2014di}}
\begin{equation}\label{eq:Hasy:12}
	A^+ := \mathbb{J}_{2k}^T A^T \mathbb{J}_{2n}.
\end{equation}
\end{definition}
We point out the properties of the symplectic inverse and refer the reader to \cite{Peng:2014di} for detailed proof.
\begin{lemma} \label{lemma:Hasy:1}
Let $A\in \mathbb{R}^{2n\times 2k}$ be a symplectic matrix and $A^+$ its symplectic inverse as defined in (\ref{eq:Hasy:12}). Then ${(A^+)}^T$ is a symplectic matrix and $A^+A = I_{2k}$.
\end{lemma}

{\edit A straight-forward calculation verifies} that $AA^+$ is idempotent, i.e., a symplectic projection onto the column span of $A$.

{\edit It is natural to expect a numerical integrator that solves (\ref{eq:Hasy:7}) to also satisfy the conservation law in Theorem \ref{theorem:Hasy:1}. Common numerical integrators e.g., Runge-Kutta methods, do not generally preserve the Hamiltonian which results in a qualitative wrong behavior of the solution \cite{Hairer:1250576}. Symplectic integrators are a class of numerical integrators for Hamiltonian systems that preserve the symplectic structure and ensure stability in long-time integration.} The Str\"omer-Verlet time stepping scheme is an example of symplectic integrators and is given by
\begin{equation} \label{eq:Hasy:13}
\begin{aligned}
	q_{n+1/2} &= q_n + \frac{\Delta t}{2} \nabla_pH(q_{n+1/2},p_n), \\
	p_{n+1} &= p_n - \frac{\Delta t}{2} \left( \nabla_qH(q_{n+1/2},p_n) + \nabla_qH(q_{n+1/2},p_{n+1}) \right),\\
	q_{n+1} &= q_{n+1/2} + \frac{\Delta t}{2} \nabla_pH(q_{n+1/2},p_{n+1}),
\end{aligned}
\end{equation}
and
\begin{equation} \label{eq:Hasy:14}
\begin{aligned}
	p_{n+1/2} &= p_n - \frac{\Delta t}{2} \nabla_qH(q_{n},p_{n+1/2}), \\
	q_{n+1} &= q_n + \frac{\Delta t}{2} \left( \nabla_pH(q_{n},p_{n+1/2}) + \nabla_pH(q_{n+1},p_{n+1/2}) \right),\\
	p_{n+1} &= p_{n+1/2} - \frac{\Delta t}{2} \nabla_qH(q_{n+1},p_{n+1/2}).
\end{aligned}
\end{equation}
For a general Hamiltonian system, the Str\"omer-Verlet scheme is implicit. However, for separable Hamiltonians, i.e. $H(q,p) = K(p) + U(q)$, this {\edit scheme becomes} explicit. We refer the reader to \cite{Hairer:1250576} for more information about the construction and applications of symplectic and geometric numerical integrators.

\section{Symplectic Model Reduction} \label{chap:SyMo:1}
We now {\edit discuss} how to modify reduced order modeling to ensure that {\edit the resulting scheme preserves} the symplectic structure of the Hamiltonian system.

Consider a Hamiltonian system (\ref{eq:Hasy:7}) on a $2n$-dimensional symplectic vector space $\edit (V,\Omega)$. Suppose that the solution manifold $\mathcal M_H$ is well approximated by a low dimensional symplectic subspace $\edit (W,\Omega)$ of dimension $2k$ $(k\ll n)$. We can {\edit then} construct a symplectic basis $A$ for $\edit W$ and approximate the solution to (\ref{eq:Hasy:7}) as
\begin{equation} \label{eq:SyMo:1}
	\mathbf z \approx A\mathbf y.
\end{equation}
Substituting this into (\ref{eq:Hasy:7}) we obtain
\begin{equation} \label{eq:SyMo:2}
	A\mathbf y = \mathbb{J}_{2n} \nabla_{\mathbf z} H(A \mathbf y). 
\end{equation}
Multiplying both sides with the symplectic inverse of $A$ and using the chain rule we have
\begin{equation} \label{eq:SyMo:3}
	\mathbf y = A^+ \mathbb J_{2n} (A^+)^T \nabla_{\mathbf y} H(A\mathbf y).
\end{equation}
Since $A$ is a symplectic basis, Lemma \ref{lemma:Hasy:1} ensures that $(A^+)^T$ is a symplectic matrix i.e., $A^+ \mathbb J_{2n} (A^+)^T = \mathbb{J}_{2k}$. By defining the reduced Hamiltonian $\tilde H:\mathbb R^{2k} \to \mathbb R$ as $\tilde H (y) = H(Ay)$ we obtain the reduced system
\begin{equation} \label{eq:SyMo:4}
\left\{
\begin{aligned}
	 \frac{d}{dt} \mathbf y &= \mathbb J_{2k} \nabla_{\mathbf y} \tilde H(\mathbf y), \\
	 \mathbf y_0 &= A^+ \mathbf z_0.
\end{aligned}
\right.
\end{equation}
The system obtained from the Petrov-Galerkin projection in (\ref{eq:MoOr:5}) is not a Hamiltonian system and does not guarantee conservation of the symplectic structure. On the other hand, we observe that the reduced system in (\ref{eq:SyMo:4}) is of the form (\ref{eq:Hasy:7}) and, hence, is a Hamiltonian system, i.e. the symplectic structure will be conserved along integral curves of (\ref{eq:SyMo:4}). Note that the original and the reduced systems are {\edit endowed with} different Hamiltonians. In the next proposition we show that the error in the Hamiltonian is constant in time.

\begin{proposition}
Let $\mathbf{z} (t)$ be the solution of (\ref{eq:Hasy:7}) at time $t$. Further suppose that $\tilde{\mathbf{z}} (t)$ is the approximate solution of the reduced system (\ref{eq:SyMo:4}) in the original coordinate system. Then the error in the Hamiltonian defined by
\begin{equation} \label{eq:SyMo:5}
	\Delta H(t)  = |H(\mathbf z(t)) - H(\tilde{\mathbf z}(t))|,
\end{equation}
is constant for all $t\in \mathbb R$.
\end{proposition}

\begin{proof}
	Let $\phi_t$ and $\psi_t$ be the Hamiltonian flow of the original and the reduced system respectively. By definition $\mathbf z(t) = \phi_t(\mathbf z_0)$ and $\mathbf y(t) = \psi_t(\mathbf y_0)$. Using the definition of the reduced Hamiltonian and Theorem \ref{theorem:Hasy:1} we have
\begin{equation} \label{eq:SyMo:6}
\begin{aligned}
	H(\tilde{\mathbf{z}} (t)) = H( A\mathbf y (t) ) = \tilde H(\mathbf y (t)) = \tilde H(\psi_t(\mathbf y_0)) = \tilde H(\mathbf y_0) = \tilde H(A^+ \mathbf z_0) = H(AA^+\mathbf z_0).
\end{aligned}
\end{equation}
The error in the Hamiltonian can then be written in terms of $\mathbf z_0$ and the symplectic basis $A$ as
\begin{equation} \label{eq:SyMo:7}
	\Delta H(t) = |H(\mathbf z_0) - H(AA^+\mathbf z_0)|
\end{equation}
\end{proof}

{\edit 
The following theorems provide a strong indication of the stability of the reduced system. 

\begin{definition} \label{definition:SyMo:1} \cite{bhatia2002stability}
Consider a dynamical system of the form $\dot{\mathbf z} = \mathbf f(\mathbf z)$ and suppose that $\mathbf z_e$ is an equilibrium point for the system so that $\mathbf f(\mathbf z_e) = 0$. $\mathbf z_e$ is called nonlinearly stable or Lyapunov stable if, for any $\epsilon > 0$, we can find $\delta > 0$ such that for any trajectory $\phi_t$, if $\| \phi_0 - \mathbf z_e \|_2 \leq \delta$, then for all $0 \leq t < \infty$, we have $\| \phi_t - \mathbf z_e \|_2 < \epsilon$, where $\| \cdot \|_2$ is the Euclidean norm.
\end{definition}	
The following proposition, also known as Dirichlet's theorem \cite{bhatia2002stability}, states the sufficient condition for an equilibrium point to be Lyapunov stable. We refer the reader to \cite{bhatia2002stability} for the proof.
\begin{proposition} \label{proposition:SyMo:1} \cite{bhatia2002stability}
An equilibrium point $\mathbf z_e$ is Lyapunov stable if there exists a scalar function $W : \mathbb R^{n} \to  \mathbb R$ such that $\nabla W(\mathbf z_e) = 0$, $\nabla^2 W(\mathbf z_e)$ is positive definite, and that for any trajectory $\phi_t$ defined in the neighborhood of $\mathbf z_e$, we have $\frac{d}{dt} W(\phi_t) \leq 0$. Here $\nabla^2W$ is the Hessian matrix of $W$.
\end{proposition}
The scalar function $W$ is referred to as the \emph{Lyapunov function}. In the context of the Hamiltonian systems, a suitable candidate for the Lyapunov function is the Hamiltonian function $H$. The following theorem shows that when $H$ (or $-H$) is a Lyapunov function, then the equilibrium points of the original and the reduced system are Lyapunov stable \cite{abraham1978foundations}. 
\begin{theorem} \label{theorem:SyMo:1}
Consider a Hamiltonian system of the form (\ref{eq:Hasy:7}) together with the reduced system (\ref{eq:SyMo:4}). Suppose $\mathbf z_e$ is an equilibrium point for (\ref{eq:Hasy:7}) and that $\mathbf y_e = A^+\mathbf z_e$. If $H$ (or $-H$) is a Lyapunov function satisfying Proposition \ref{proposition:SyMo:1}, then $\mathbf z_e$ and $\mathbf y_e$ are Lyapunov stable equilibrium points for (\ref{eq:Hasy:7}) and (\ref{eq:SyMo:4}), respectively. 
\end{theorem}
\begin{proof}
	It is a direct consequence of Proposition \ref{proposition:SyMo:1} that $\mathbf z_e$ is a local minimum or maximum of (\ref{eq:Hasy:7}) and also a Lyapunov stable point. It can be easily checked that if $\mathbf z_e$ is a local minimum of $H$ then $\mathbf y_e$ is a local minimum for $\tilde H$ and an equilibrium point for (\ref{eq:SyMo:4}). Also from the chain rule we have
\[
	\nabla^2_{\mathbf y} \tilde H = A^T \nabla^2_{\mathbf z} H A.
\]
So for any $\xi\in \mathbb R^{2k}$
\[
	\xi^T \nabla^2_{\mathbf y} \tilde H \xi = (A\xi)^T \nabla^2_{\mathbf z} H (A\xi) \geq 0.
\]
Here the last inequality is due to the positive definiteness of $H$. Therefore $\tilde H$ is also positive definite. By Proposition \ref{proposition:SyMo:1} we conclude that $\mathbf y_e$ is a Lyapunov stable point.
\end{proof}
}

While the symplectic structure is not guaranteed to be preserved in the reduced systems obtained by the Petrov-Galerkin projection, the reduced system obtained by the symplectic projection guarantees the preservation of the energy up to the error in the Hamiltonian (\ref{eq:SyMo:5}). In the next section we discuss  different methods for obtaining a symplectic basis.

\subsection{Proper Symplectic Decomposition (PSD)} \label{chap:SyMo.PrSy:1}

Similar to Section \ref{chap:MoOr.PrOr:1} we gather snapshots $\mathbf z_i = [q_i^T , p_i^T]^T$ in the snapshot matrix $S$. Suppose that a symplectic basis $A$ of size $2n\times2k$ and its symplectic inverse $A^+$ is provided. {\edit The Proper Symplectic Decomposition} requires that the error of the symplectic projection onto the symplectic subspace {\edit be minimized}. Hence, the PSD symplectic basis of size $2k$ is the solution to the optimization problem

\begin{equation} \label{eq:SyMo:8}
\begin{aligned}
& \underset{V\in \mathbb R^{2n\times 2k}}{\text{minimize}}
& & \| S - AA^+S\|_F \\
& \text{subject to}
& & A^T \mathbb{J}_{2n}A = \mathbb{J}_{2k}
\end{aligned}
\end{equation}
Compared to POD, in (\ref{eq:SyMo:8}) the orthogonal projection is replaced with a symplectic projection $AA^+$. At first, the minimization looks similar to the one obtained by POD. {\edit However, it is well known that symplectic bases are not generally orthogonal, and therefore not norm bounded. This means that numerical errors may become dominant in the symplectic projection \cite{Karow:2006cf} which makes the minimization (\ref{eq:SyMo:8}) a harder problem than (\ref{eq:MoOr:6}).}
	
As the optimization problem (\ref{eq:SyMo:8}) is nonlinear, the direct solution is usually expensive. A simplified version of the optimization (\ref{eq:SyMo:8}) can be found in \cite{Peng:2014di}, but there is no guarantee that the method provides a near optimal basis. 

{\edit Finding eigen-spaces of Hamiltonian and symplectic matrices is studied in the context of optimal control problems \cite{Benner:2000ww,Benner:1997ef,Watkins:2004kz,BunseGerstner:1986dg} and model reduction of Riccati equations \cite{Benner:1997ef}, where also an SVD-like decomposition for Hamiltonian and symplectic matrices has been proposed \cite{Xu:2003kx}. However, the computation of a large snapshot matrix and use of the mentioned methods to compute its eigen-spaces, is usually computationally demanding. Also, these methods generally do not guarantee the construction of a well-conditioned symplectic basis.
	
The greedy approach presented in Section \ref{Chap:Symo.PrSy:3} is an iterative method for construction of a symplectic basis. It avoids the evaluation of the full snapshot matrix, hence substantially reduces the computational cost in the offline stage of the symplectic model reduction. Also, by construction, it yields an orthosymplectic basis and therefore a well-conditioned basis.

In Section \ref{chap:SyMo.PrSy:2} we briefly outline non-direct methods for finding solutions to (\ref{eq:SyMo:8}), proposed by \cite{Peng:2014di}, and assuming a specific structure for $A$. In Section \ref{Chap:Symo.PrSy:3} we introduce a greedy approach for the symplectic basis generation.}

\subsubsection{SVD Based Methods for Symplectic Basis Generation} \label{chap:SyMo.PrSy:2}

\paragraph{\bf Cotangent lift} Suppose that $A$ is of the form
\begin{equation} \label{eq:SyMo:9}
	A = 
	\begin{pmatrix}
		\Phi & 0 \\
		0 & \Phi
	\end{pmatrix},
\end{equation}
where $\Phi \in \mathbb{R}^{n\times k}$ is an orthonormal matrix. It is easy to check that $A$ is a symplectic matrix, i.e., $A^T \mathbb J_{2n} A = \mathbb J_{2k}$. The construction of $A$ suggests that the range of $\Phi$ should cover both the potential and the momentum spaces. Hence, we can construct $A$ by forming the combined snapshot matrix
\begin{equation} \label{eq:SyMo:10}
	S_{\text{combined}} = [q_1,\dots,q_n,p_1,\dots,p_n], \qquad \mathbf z_i = (q_i^T,p_i^T)^T,
\end{equation}
and define $\Phi=[u_1,\dots,u_k]$, where $u_i$ is the $i$-th left singular vector of $S_{\text{combined}}$. It is shown in \cite{Peng:2014di} that among all symplectic bases of the form (\ref{eq:SyMo:9}) cotangent lift minimizes the projection error.

\paragraph{\bf Complex SVD} Suppose instead that $A$ takes the form \cite{Peng:2014di}
\begin{equation} \label{eq:SyMo:11}
	A = 
	\begin{pmatrix}
		\Phi & -\Psi \\
		\Psi & \Phi
	\end{pmatrix},
\end{equation}
while $\Phi$ and $\Psi$ are real matrices of size $n\times k$ satisfying conditions
\begin{equation} \label{eq:SyMo:12}
\Phi^T \Phi + \Psi^T \Psi = I_k,\quad \Phi^T \Psi = \Psi^T \Phi.
\end{equation}
It can be checked that $A$ forms a symplectic matrix. To construct $A$ we first define the complex snapshot matrix
\begin{equation} \label{eq:SyMo:13}
	S_{\text{complex}} = [ q_1 + i p_1, \dots , q_N + i p_N ].
\end{equation}
Each left singular vector of $S_{\text{complex}}$ now takes the form $u_m = r_m + i s_m$. We define
\begin{equation} \label{eq:SyMo:14}
	 \Phi = [r_1,\dots, r_k], \quad \Psi = [s_1,\dots, s_k].
\end{equation}
One can easily check that (\ref{eq:SyMo:12}) is satisfied {\edit since the matrix of singular vectors is unitary}. It is shown in \cite{Peng:2014di} that among all symplectic bases of the form (\ref{eq:SyMo:11}) the complex SVD minimizes the projection error.

\subsubsection{The Greedy Approach to Symplectic Basis Generation} \label{Chap:Symo.PrSy:3} Greedy generation of the reduced basis is an iterative procedure which, in each iteration, adds the two best possible basis vectors to the symplectic basis to enhance overall accuracy. In contrast to the cotangent lift and the complex SVD methods, the greedy approach does not require the symplectic basis to have a specific structure. This typically results in a more compact basis and/or more accurate reduced systems. For parametric problems, the greedy approach only requires one numerical solution to be computed per iteration hence saving substantial computational cost in the offline stage. 

{\edit The orthonormalization step is an essential step in most greedy approaches for basis generation in the context of model reduction \cite{Anonymous:2016wl,Quarteroni:2016wi}. However common orthonormalization processes, e.g. the QR method, destroy the symplectic structure of the original system \cite{BunseGerstner:1986dg}. Here we use a variation of the QR method known as the SR \cite{Salam2014} method which is based on the symplectic Gram-Schmidt method and yields a symplectic basis. 
}

{\edit
As discussed in Section \ref{chap:Hasy:1}, any finite dimensional symplectic linear vector space has a symplectic basis that satisfies conditions (\ref{eq:Hasy:4}). Further, Theorem \ref{theorem:Hasy:1.6} provides an iterative process for constructing an orthosymplectic basis \cite{Matsuo:2014wl,Salam2014}. To briefly describe the SR method, suppose that an orthosymplectic basis
\begin{equation} \label{eq:SyMo:14.1}
	A_{2k}=\{ e_1 , \dots , e_k \} \cup \{ \mathbb J_{2n}^T e_1 , \dots , \mathbb J_{2n}^T e_k \},
\end{equation}
and a vector $z\not \in \text{span}(A_{2k})$ is provided. We aim to symplectically orthogonalize ($\mathbb J_{2n}$-orthogonalize) $z$ with respect to $A_{2k}$ and seek $\alpha_1,\dots,\alpha_k,\beta_1,\dots,\beta_k \in \mathbb R$ such that
\begin{equation} \label{eq:SyMo:14.2}
	\Omega\left(z + \sum_{i=1}^k \alpha_i e_i + \sum_{i=1}^k \beta_i \mathbb J_{2n}^Te_i , \sum_{i=1}^k \bar{\alpha}_i e_i + \sum_{i=1}^k \bar{\beta}_i \mathbb J_{2n}^Te_i \right) = 0,
\end{equation}
for all possible $\bar{\alpha}_1,\dots,\bar{\alpha}_k,\bar{\beta}_1,\dots,\bar{\beta}_k \in \mathbb R$. It is easily seen that the unique solution is 
\begin{equation} \label{eq:SyMo:14.3}
	\alpha_i = - \Omega(z,\mathbb J_{2n}^Te_i), \quad \beta_i = \Omega(z,e_i),
\end{equation}
for $i=1,\dots,k$. Now define the modified vectors as
\begin{equation} \label{eq:SyMo:14.4}
	\tilde z = z - \sum_{i=1}^k \Omega(z,\mathbb J_{2n}^Te_i) e_i + \sum_{i=1}^k \Omega(z,e_i) \mathbb J_{2n}^Te_i.
\end{equation}
If we introduce $e_{k+1} = \tilde z / \| \tilde z \|_2$, it is easily checked that $e_{k+1}$ is also orthogonal to $A_{2k}$ with respect to the classical inner product. Therefore span$\{e_1,\dots,e_{k+1}\}$ forms a Lagrangian subspace and according to Theorem \ref{theorem:Hasy:1.6} the basis $A_{2k+2}= A_{2k}\cup \{ e_{k+1} , \mathbb J_{2n}^T e_{k+1} \}$ forms an orthosymplectic basis.

Note that the $SR$ method can be replaced with backward stable routines such as the isotropic Arnoldi or the isotropic Lanczos methods \cite{Mehrmann:2000dv}.
}

The key element of the greedy algorithm is the availability of an error function which evaluates the error associated with the model reduction \cite{Anonymous:2016wl}. In the framework of symplectic model reduction, one possible candidate is the error in the Hamiltonian (\ref{eq:SyMo:5}). Correctly approximating symplectic systems relies on preservation of the Hamiltonian, hence the error in the Hamiltonian {\edit arises as a} a natural choice. Moreover, since the error in the Hamiltonian depends on the initial condition and the reduced symplectic basis, evaluation of the error does not require the time integration of the full system. 

Suppose that a $2k$-dimensional {\edit orthosymplectic basis (\ref{eq:SyMo:14.1})} is generated at the $k$-th step of the greedy method and we seek to enrich it by two additional vectors. Using the error in the Hamiltonian (\ref{eq:SyMo:7}) we search the parameter space to identify the value that maximizes the error in the Hamiltonian
\begin{equation} \label{eq:SyMo:14.5}
	\omega_{k+1} := \underset{\omega\in \Gamma}{\text{argmax }}\Delta H(\omega).
\end{equation}
The goal is to approximate the Hamiltonian function as well as possible. 

We then propagate (\ref{eq:Hasy:7}) in time to produce trajectory snapshots 
\begin{equation}
	S=\{ \mathbf z(t_i,\omega_{k+1}) | i = 1,\dots,M \}.
\end{equation} 
The next basis vector is the snapshot that maximises the projection error (\ref{eq:SyMo:8})
{\edit 
\begin{equation} \label{eq:SyMo:14.6}
	z := \underset{s\in S}{\text{argmax }} \| s - A_{2k}{A_{2k}}^+s \|.
\end{equation}
}
Finally, we update the basis as
{\edit
\begin{equation} \label{eq:SyMo:14.7}
	e_{k+1} = \tilde z, \quad A_{2k+1} = A_{2k}\cup \{ e_{k+1} , \mathbb J_{2n}^Te_{k+1} \},
\end{equation}
}
where $\tilde z$ is the vector obtained {\edit after} applying the symplectic Gram-Schmidt process to $z$. 

Since the maximization over the entire parameter space $\Gamma$ is impossible, we discretize the parameter set into a grid with $N$ points: $\Gamma_N = \{ \omega_1,\dots,\omega_N\}$. However, since the selection of parameters only require the evaluation of the error in the Hamiltonian and not time integration of the original system, then $\Gamma_N$ can be chosen {\edit to be} very rich.

We summarize the greedy algorithm for the generation of a symplectic basis in Algorithm \ref{alg:SyMo:3}.

\begin{algorithm} 
\caption{The greedy algorithm for generation of a symplectic basis} \label{alg:SyMo:3}
{\bf Input:} Tolerated loss in the Hamiltonian $\delta$, parameter set $\Gamma_N = \{\omega_1,\dots,\omega_N\}$, initial condition $\mathbf z_0(\omega)$
\begin{enumerate}
\item $\omega^* \leftarrow \omega_1$
\item $e_1 \leftarrow \mathbf z_0(\omega^*)$
\item $A \leftarrow [e_1,\mathbb J^T_{2n}e_1]$
\item $k \leftarrow 1$
\item \textbf{while} $\Delta H(\omega) > \delta$ for all $\omega \in \Gamma_N$
\item \hspace{0.5cm} $\omega^* \leftarrow$ $\underset{\omega\in \Gamma_N}{\text{argmax }}\Delta H(\omega)$
\item \hspace{0.5cm} Compute trajectory snapshots $S=\{ \mathbf z(t_i,\omega^*) | i = 1,\dots,M \}$
\item \hspace{0.5cm} $\mathbf z^* \leftarrow$ $\underset{s\in S}{\text{argmax }} \| s - AA^+s \|$
\item \hspace{0.5cm} Apply symplectic Gram-Schmidt on $\mathbf z^*$
\item \hspace{0.5cm} $e_{k+1} \leftarrow \mathbf z^*/ \| \mathbf  z^*\|$
\item \hspace{0.5cm} $A \leftarrow [e_1,\dots ,e_{k+1} , \mathbb J^T_{2n}e_1,\dots,\mathbb J^T_{2n}e_{k+1}]$
\item \hspace{0.5cm} $k \leftarrow k+1$
\item \textbf{end while}
\end{enumerate}
\vspace{0.5cm}
{\bf Output:} Symplectic basis $A$.
\end{algorithm}


\subsubsection{Convergence of the Greedy Method} \label{chap:SyMo.PrSy:3}

To show convergence of the greedy method we {\edit consider} a slightly different version based on the projection error. The error in the Hamiltonian is then introduced as a cheap surrogate to the projection error to accelerate the parameter selection.

Suppose that we are given a compact subset $S$ of $\mathbb R^{2n}$. Our intention is to find a set of vectors $A=\{e_1,\dots,e_k,f_1,\dots,f_k\}$ such that $A$ forms {\edit an orthosymplectic} basis and any $s\in S$ is well approximated by elements of the subspace span$(A)$. The modified greedy method for generating basis vectors $e_i$ and $f_i$ is as follows. In the initial step we pick $e_1$ such that $\edit \|e_1\|_2 = \max_{s\in S} \|s\|_2$. Then define $f_1 = \mathbb{J}_{2n}^T e_1$. It is easy to check that the span of $A_2 = \{e_1,f_1\}$ is {\edit orthosymplectic}, so $A_2$ is the first subspace that approximates elements of $S$. In the $k$-th step of the greedy method, suppose we have a basis $A_{2k} = \{ e_1,\dots, e_k , f_1,\dots ,f_k \}$. We define $P_{2k}$ to be a symplectic projection operator that projects elements of $S$ onto span$(A_{2k})$ and define
\begin{equation} \label{eq:new1}
	\sigma_{2k}(s) := \|s-P_{2k}(s)\|_2,
\end{equation}
as the projection error. Moreover we denote by $\sigma_{2k}$ the maximum approximation error of $S$ using elements in span$(A_{2k})$ as
\begin{equation} \label{eq:new2}
	\sigma_{2k} := \max_{s\in S} \sigma_{2k}(s).
\end{equation}
The next set of basis vectors in the greedy selection are
\begin{equation} \label{eq:new3}
	e_{k+1} := \underset{s\in S}{\text{argmax }}\sigma_{2k}(s), \quad f_{k+1} := \mathbb{J}_{2n}^T e_{k+1}.
\end{equation}
We emphasisze that the sequence of basis vectors generated by the greedy is generally not unique. 

To estimate the quality of the reduced subspace, it is natural to compare it with the best possible $2k$-dimensional subspace in the sense of the minimum projection (not necessary symplectic) error. For this we introduce the Kolmogorov $n$-width \cite{Kolmogoroff:1936fj,Pinkus:1985vy}.

\begin{definition}
Let $S$ be a subset of $\mathbb R^{m}$ and $Y_n$, $n\leq m$, be a general $n$-dimensional subspace of $\mathbb R^{m}$. The angle between $S$ and $Y_n$ is given by
\begin{equation} \label{eq:new4}
	E(S,Y_n) := \sup_{s\in S} \inf_{y\in Y_n} \|s-y\|_2.
\end{equation}
The Kolmogorov $n$-width of $S$ in $\mathbb R^m$ is given by
\begin{equation} \label{eq:new5}
	d_{n}(S,\mathbb{R}^m) := \inf_{Y_n} E(S,Y_n) = \inf_{Y_n} \sup_{s\in S} \inf_{y\in Y_n} \|s-y\|_2
\end{equation}
\end{definition}

For a given subspace $Y_n$, the angle between $S$ and $Y_n$ measures the worst possible projection error of elements in $S$ onto $Y_n$. Hence the Kolmogorov $n$-width quantifies how well $S$ can be approximated by {\edit an} $n$-dimensional subspace. 

We seek to show that the decay of $\sigma_{2k}$, obtained by the greedy algorithm, has the same rate as of $d_{2k}(S)$, i.e., the greedy method provides the best possible accuracy attained by a $2k$-dimensional subspace.

We start by $\edit \mathbb J_{2n}$-orthogonalizing the vectors provided by the greedy algorithm as
\begin{equation} \label{eq:new6}
\begin{aligned}
	& \xi_1 = e_i, & \bar{\xi}_1 = \mathbb{J}_{2n}^T \xi_1, &\\
	& \xi_i = e_i - P_{2(i-1)} (e_i), & \bar{\xi}_i = \mathbb{J}_{2n}^T, \xi_i &\quad i = 2,3,\dots
\end{aligned}
\end{equation}
The projection of a vector $s\in S$ onto span$(A_{2k})$ can be written using the symplectic basis as
\begin{equation} \label{eq:new7}
	P_{2k}(s) = \sum_{i=1}^k \left( \alpha_i(s) \xi_i + \bar{\alpha}_i(s) \bar{\xi}_i \right),
\end{equation}
where $\alpha_i(s)$ and $\bar{\alpha}_i(s)$ for $i=1,\dots,k$ are the expansion coefficients
\begin{equation} \label{eq:new8}
	\alpha_i(s) = - \frac{\Omega(\bar{\xi}_i,s)}{\Omega(\xi_i,\bar{\xi}_i)}, \quad \bar{\alpha}_i(s) = \frac{\Omega(\xi_i,s)}{\Omega(\xi_i,\bar{\xi}_i)},
\end{equation}
for any $s\in S$. Since $\bar{\xi}_i$ is {\edit $\mathbb J_{2n}$-orthogonal} to the span$(A_{2(k-1)})$ we have
\begin{equation} \label{eq:new9}
\begin{aligned}
	|\alpha_i(s)| = \frac{|\Omega(\bar{\xi}_i,s)|}{|\Omega(\xi_i,\bar{\xi}_i)|} = \frac{|\Omega( \bar{\xi}_i, s - P_{2(k-1)}(s))|}{|\Omega(\xi_i,\bar{\xi}_i)|}  &\leq \frac{\|\bar{\xi}_i\|_2 \| s - P_{2(k-1)}(s) \|_2}{ \|\xi_i\|_2 \|\bar{\xi}_i\|_2 } \\
	&= \frac{\| s - P_{2(k-1)}(s) \|_2}{\| e_i - P_{2(k-1)}(e_i) \|_2} \leq 1.
\end{aligned}
\end{equation}
Here, we use the fact that $\edit |\Omega(\xi_i,\bar{\xi}_i)| = \| \xi_i \|^2_2 = \|\bar{\xi}_i\|^2_2$ with the last inequality following from the greedy algorithm which maximizes $e_i$. Similarly we deduce that $|\bar{\alpha}_i(s)|\leq 1$.

We write
\begin{equation} \label{eq:new10}
\begin{aligned}
	\xi_j = \sum_{i=1}^j \left( \mu_i^j e_i + \gamma_i^j f_i \right), \quad \bar{\xi}_j = \sum_{i=1}^j \left( \lambda_i^j e_i + \eta_i^j f_i,  \right), \quad j=1,2,\dots
\end{aligned}
\end{equation}
with
\begin{equation} \label{eq:new11}
\begin{aligned}
	&\mu^j_j = 1, \quad \gamma^j_j = 0, \\
	&\mu_i^j = \sum_{l=i}^{j-1}\left( - \alpha_l(f_j) \mu_i^l  + \bar{\alpha}_l(f_j) \gamma_i^l \right), \quad\gamma_i^j = \sum_{l=i}^{j-1}\left( - \alpha_l(f_j) \gamma_i^l  + \bar{\alpha}_l(f_j) \mu_i^l \right), \\
	&\lambda^j_i = - \gamma ^j_i, \quad \eta^j_i = \mu^j_i,
\end{aligned}
\end{equation}
for $j=2,3,\dots$. By induction and using the bound in (\ref{eq:new9}) we deduce that
\begin{equation} \label{eq:new12}
	\mu^j_i,\gamma^j_i,\lambda^j_i,\eta^j_i \leq 3^{j-i}, \quad \text{for } j\geq i.
\end{equation}
Now let $2k$ be the dimension of the desired reduced space. Looking at the definition of Kolmogorov $n$-width we observe that for any $\theta > 1$ we can find a subspace $Y_{2k}$ such that $E(S,Y_{2k}) \leq \theta d_{2k}(S,\mathbb R^n)$. Hence we can find vectors $v_1,\dots,v_k,u_1,\dots,u_k\in Y_{2k}$ such that
\begin{equation} \label{eq:new13}
\begin{aligned}
	& \|e_i - v_i\|_2 \leq \theta d_{2k}(S,\mathbb R^n), \\
	& \|f_i - u_i\|_2 \leq \theta d_{2k}(S,\mathbb R^n).
\end{aligned}
\end{equation}
Now we construct a set of $2(k+1)$ new vectors
\begin{equation} \label{eq:new14}
\begin{aligned}
	& \zeta_j = \sum_{i=1}^{k+1} \mu_i^j v_i + \gamma^j_i u_i,\quad \bar{\zeta}_j = \sum_{i=1}^{k+1} \lambda_i^j v_i + \eta^j_i u_i.
\end{aligned}
\end{equation}
for $j = 1,\dots,k+1$. Note that since $u_i$ and $v_i$ belong to $Y_{2k}$ so does their linear combination including all $\zeta_j$ and $\bar{\zeta}_j$. We can use the inequality (\ref{eq:new12}) to write
\begin{equation}
	\| \xi_i - \zeta_i \|_2 \leq 3^i \theta d_{2k}(S,\mathbb R^n),\quad \| \bar{\xi}_i - \bar{\zeta}_i \|_2 \leq 3^i \theta d_{2k}(S,\mathbb R^n).
\end{equation}
Moreover since $Y_{2k}$ is of dimension $2k$ we find $\kappa_i$, $i=1,\dots,2(k+1)$ such that
\begin{equation} \label{eq:new15}
	\sum_{i=1}^{2(k+1)} \kappa_i^2 = 1, \quad\sum_{i=1}^{k+1} \kappa_i \zeta_i + \sum_{i=1}^{k+1} \kappa_{i+k+1} \bar{\zeta}_i = 0.
\end{equation}
We have
\begin{equation} \label{eq:new17}
\begin{aligned}
	\left\| \sum_{i=1}^{k+1} \kappa_i \xi_i + \sum_{i=1}^{k+1} \kappa_{i+k+1} \bar{\xi}_i \right\|_2 &= \left\| \sum_{i=1}^{k+1} \kappa_i (\xi_i - \zeta_i) + \sum_{i=1}^{k+1} \kappa_{i+k+1} (\bar{\xi}_i-\bar{\zeta}_i) \right\|_2 \\
	&\leq 2\cdot 3^{k+1} \sqrt{2(k+1)} \theta d_{2k}(S,\mathbb R^n).
\end{aligned}
\end{equation}
We know there exists $1 \leq j\leq 2k+2$ such that $\kappa_j > 1/\sqrt{2(k+1)}$. Without loss of generality let us assume that $j\leq k+1$. This yields
\begin{equation} \label{eq:new18}
	\left\| \xi_j +  \kappa_j^{-1} \sum_{i=1,i\neq j}^{k+1} \kappa_i \xi_i + \kappa_j^{-1}\sum_{i=1}^{k+1} \kappa_{i+k+1} \bar{\xi}_i \right\|_2 \leq 4\cdot 3^{k+1} (k+1) \theta d_{2k}(S,\mathbb R^n).
\end{equation}
Define $c = \kappa_j^{-1} \sum_{i=1,i\neq j}^{k+1} \kappa_i \xi_i + \kappa_j^{-1}\sum_{i=1}^{k+1} \kappa_{i+k+1} \bar{\xi}_i$. Using that $\mathbb{J}_{2n}^T c$ is {\edit $\mathbb J_{2n}$-orthogonal} to $\xi_j$ we recover
\begin{equation} \label{eq:new19}
\begin{aligned}
	\| \xi_j \|_2 &\leq \| \xi_j \|_2 + \| c \|_2 = \Omega(\xi_j,\mathbb{J}_{2n}^T \xi_j) + \Omega(c,\mathbb{J}_{2n}^T c) \\
	       &= \Omega(\xi_j,\mathbb{J}_{2n}^T \xi_j) + \Omega(c,\mathbb{J}_{2n}^T c) + \Omega(\xi_j,\mathbb{J}_{2n}^T c) + \Omega(c,\mathbb{J}_{2n}^T \xi_j) \\
       	&= \Omega(\xi_j + c, \mathbb{J}^T_{2n} (\xi_j + c)) = \| \xi_j + c \|_2	
\end{aligned}
\end{equation}
Combining this with (\ref{eq:new18}) yields
\begin{equation} \label{eq:new20}
	\| \xi_j \|_2 \leq 4\cdot 3^{k+1} (k+1) \theta d_{2k}(S,\mathbb R^n).
\end{equation}
Finally using the definition of $\xi_j$ for all $s\in S$ we have
\begin{equation} \label{eq:new21}
	\| s - P_{2(j-1)}(s) \|_2 \leq \| f_j - P_{2(j-1)}(f_j) \|_2 = \|\xi_j \|_2 \leq 4\cdot 3^{k+1} (k+1) \theta d_{2k}(S,\mathbb R^n)
\end{equation}
Hence, for any given $\lambda > 1$
\begin{equation} \label{eq:new22}
	\| s - P_{2k}(s) \|_2 \leq \| s - P_{2(j-1)}(s) \|_2 \leq 4\cdot 3^{k+1} (k+1) \theta d_{2k}(S,\mathbb R^n).
\end{equation}
This establishes the following theorem.
\begin{theorem} \label{theorem:SyMo:2}
	Let $S$ be a compact subset of $\mathbb{R}^{2n}$ with exponentially small Kolmogorov $n$-width $\edit d_{k}\leq c\exp(-\alpha k)$ with $\alpha > \log3$. Then there exists $\beta>0$ such that the symplectic subspaces $A_{2k}$ generated by the greedy algorithm provide exponential approximation properties such that
\begin{equation} \label{eq:new23}
	\| s - P_{2k}(s) \|_2 \leq C \exp(-\beta k)
\end{equation}
for all $s\in S$ and some $C>0$.
\end{theorem}


\subsection{Symplectic Discrete Empirical Interpolation Method (SDEIM)} Consider the Hamiltonian system (\ref{eq:Hasy:7}) and its reduced system (\ref{eq:SyMo:4}) equipped with a symplectic transformation $A$. One can split the Hamiltonian function $H = H_1 + H_2$ such that $\nabla H_1 = L\mathbf z$ and $\nabla H_2 = \mathbf g(\mathbf z)$, where $L$ is a constant matrix in $\edit \mathbb R^{2n\times 2n}$ and $\mathbf g$ is a nonlinear function. The reduced system takes the form
\begin{equation} \label{eq:new24}
	\frac{d}{dt} \mathbf y = \underbrace{A^+ \mathbb J_{2n} L A}_{\tilde L} \mathbf y + A^+ \mathbb J_{2n} \mathbf g(A\mathbf y)
\end{equation}
As discussed in Section \ref{chap:MoOr.DEIM:1}, the complexity of evaluating the nonlinear term still depends on $n$, the size of the original system. To overcome this computational bottleneck we use the DEIM approximation for evaluating the nonlinear function $\mathbf g$ as
\begin{equation} \label{eq:new25}
	\frac{d}{dt} \mathbf y = \tilde L \mathbf y + \underbrace{ A^+ \mathbb J_{2n} V (P^TV)^{-1} P^T \mathbf g(A\mathbf y) }_{\tilde N(\mathbf y)}
\end{equation}
For a general choice of $V$ the system (\ref{eq:new25}) is not guaranteed to be a Hamiltonian system, impacting long time accuracy and stability. However, we can guarantee that (\ref{eq:new25}) is a Hamiltonian system by choosing $V=(A^+)^T$. To {\edit see this}, we note that the system (\ref{eq:new25}) is a Hamiltonian system if and only if $\tilde N(\mathbf y) = \mathbb J_{2k} \nabla_{\mathbf y} \mathbf g(\mathbf y)$. Also we have 
\begin{equation} \label{eq:new26}
	\mathbf g(A\mathbf y) = \nabla_{\mathbf z} H_2(\mathbf z) = (A^+)^T \nabla_{\mathbf y} H_2(A \mathbf y),
\end{equation}
where the chain rule is used for the second equality. Substituting this into $\tilde N$ we obtain
\begin{equation} \label{eq:new27}
	\tilde N(\mathbf y)= A^+ \mathbb J_{2n} V (P^TV)^{-1} P^T  (A^+)^T \nabla_{\mathbf y} H_2(A \mathbf y).
\end{equation}
Taking $V = (A^+)^T$ yields
\begin{equation} \label{eq:new28}
	\tilde N(\mathbf y) = A^+ \mathbb J_{2n}(A^+)^T \nabla_{\mathbf y} H_2(A \mathbf y) = \mathbb J_{2k} \nabla_{\mathbf y} H_2(A \mathbf y),
\end{equation}
since $(A^+)^T$ is a symplectic matrix. Hence, $V = (A^+)^T$ is a sufficient condition for (\ref{eq:new25}) to {\edit be Hamiltonian}. 

Regarding the construction of the projection space, suppose that we have already constructed a symplectic basis $A=\{ e_1,\dots , e_k,f_1,\dots f_k \}$ using the greedy algorithm. Note that $(A^+)^T$ is a symplectic basis and $(A^+)^+=A$. Thus, we can move between these two symplectic bases by simply using the transpose operator and the symplectic inverse operator. Let $S_{\mathbf g} = \{ \mathbf g (\mathbf x(t_i,\omega_j)) \}$ with $i = 1,\dots,M$ and $ j = 1 ,\dots,N$ be the nonlinear snapshots that were gathered in the greedy algorithm. We then form $(A^+)^T = \{ e'_1,\dots, e'_k,f'_1,\dots,f'_k\}$ and use a greedy approach to add new basis vectors to $(A^+)^T$. At the $i$-th iteration of the symplectic DEIM, we use $(A^+)^T$ to approximate elements in $S_{\mathbf g}$ and choose the vector that maximizes the error as the next basis vector 
\begin{equation}
	s^* := \underset{s \in S_{\mathbf g}}{\text{argmax }}\| s - (A^+)^T A^+ s \|_2.	
\end{equation}
After applying the symplectic Gram-Schmidt on $s^*$, we update $(A^+)^T$ as
\begin{equation}
\begin{aligned}
	e'_{k+i+1} = \frac{s^*}{\| s^* \|_2},\quad f'_{k+i+1} = \mathbb J_{2n}^T e'_{k+i+1}.
\end{aligned}
\end{equation}
Finally when $(A^+)^T$ approximates elements $S_{\mathbf g}$ with the desired accuracy, we transpose and symplectically invert $(A^+)^T$ to obtain $A$. We summarize the symplectic DEIM algorithm in Algorithm \ref{alg:SyMo:4}.

\begin{algorithm} 
\caption{Symplectic Discrete Empirical Interpolation Method} \label{alg:SyMo:4}
{\bf Input:} Symplectic basis $A=\{ e_1,\dots,e_k,f_1,\dots,f_k \}$, nonlinear snapshots $S_{\mathbf g} = \{ \mathbf g(\mathbf x(t_i,\omega_j)) \}$ and tolerance $\delta$
\begin{enumerate}
\item Compute $(A^+)^T = \{ e'_1,\dots,e'_k,f'_1,\dots,f'_k \}$
\item $i \leftarrow 1$
\item \textbf{while} max$\| s - (A^+)^T A^+s \| > \delta$ for all $s\in S_{\mathbf g}$
\item \hspace{0.5cm} $s^* \leftarrow \underset{s \in S_{\mathbf g}}{\text{argmax }}\| s - (A^+)^T A^+ s \|$
\item \hspace{0.5cm} Apply symplectic Gram-Schmidt on $s^*$
\item \hspace{0.5cm} $e'_{k+i} = s^* / \| s^* \|$
\item \hspace{0.5cm} $f'_{k+i} = \mathbb J_{2n} e'_{k+i}$
\item \hspace{0.5cm} $(A^+)^T \leftarrow [e'_1,\dots,e'_{k+i},f'_1,\dots,f'_{k+i}]$
\item \hspace{0.5cm} $i\leftarrow i+1$
\item \textbf{end while}
\item take transpose and symplectic inverse of $(A^+)^T$
\end{enumerate}
\vspace{0.5cm}
{\bf Output:} Symplectic basis $A$ that guarantees a Hamiltonian reduced system.
\end{algorithm}

When using an implicit time integration scheme we face inefficiencies when evaluating the Jacobian of nonlinear terms, as discussed in Section \ref{chap:MoOr.DEIM:1}. We recall that the key to fast approximation of the Jacobian is that the interpolating index matrix $P$, obtained in the DEIM approximation, commutes with the nonlinear function. Nonlinear terms in Hamiltonian systems often take the from
\begin{equation}
	\mathbf g (\mathbf z) = \mathbf g (\mathbf q,\mathbf p) = 
	\begin{pmatrix}
		g_1(q_1,p_1) \\
		g_2(q_2,p_2) \\
		\vdots \\
		g_{2n}(q_{n},p_{n})
	\end{pmatrix}.
\end{equation}
Thus, the interpolating index matrix, obtained by Algorithm \ref{alg:MoOr:1} does not necessarily commute with the function $\mathbf g$. To overcome this, when index $\mathfrak p_i$ with $\mathfrak p_i\leq n$ or $\mathfrak p_i>n$ is chosen in Algorithm \ref{alg:MoOr:1} we also include $\mathfrak p_i + n$ or $\mathfrak p_i-n$, respectively. {\edit Simple calculations verifies that $\mathbf g$ and $P$ then commute.}

\section{Numerical Results} \label{chap:NuRe:1}
In this section, we illustrate the performance of the greedy generation of a symplectic basis. The parametric linear wave equation is considered to compare SVD based methods with the greedy method. The nonlinear model order reduction using the combination of DIEM and the symplectic basis is then illustrated by considering the parametric nonlinear Schr\"odinger equation. {\edit Finally we discuss the numerical convergence of the greedy method introduced in Algorithm \ref{alg:SyMo:3}.}

\subsection{Parametric Linear Wave equation} \label{chap:NuRe:1.1} Consider the {\edit parametric} linear wave equation
\begin{equation} \label{eq:NuRe:1}
\left\{
\begin{aligned}
& u_{tt}(x,t,\omega) = \kappa(\omega) u_{xx}(x,t,\omega), \\
& u(x,0) = u^0(x),
\end{aligned}
\right.
\end{equation}
where $x$ belongs to a one-dimensional torus of length $L$, $\omega = (\omega_1,\dots,\omega_4)$ and
\begin{equation} \label{eq:NuRe:2}
	\kappa(\omega) = c^2\left( \sum_{l=1}^4 \frac{1}{l^2} \omega_l \right).
\end{equation}
{\edit Here $\omega_l \in [0,1]$ for $l=1,\dots,4$} and $c\in \mathbb{R}$ is a constant number. By rewriting (\ref{eq:NuRe:1}) in canonical form, using the change of variable $q = u$ and $\partial q/ \partial t= p$, we obtain the symplectic form
\begin{equation} \label{eq:NuRe:3}
\left\{
\begin{aligned}
& q_t(x,t,\omega) = p(x,t,\omega), \\
& p_t(x,t,\omega) = \kappa(\omega) q_{xx}(x,t,\omega),
\end{aligned}
\right.
\end{equation}
with the associated Hamiltonian
\begin{equation} \label{eq:NuRe:4}
	H(q,p,\omega) = \frac 1 2 \int_0^L p^2 + \kappa(\omega) q_x^2 \ dx.
\end{equation}
We discretize the torus into $N$ equidistant points and define $\Delta x = L/N$, $x_i = i\Delta x$, $q_i=q(t,x_i,\omega)$ and $p_i=p(t,x_i,\omega)$ for $i = 1, \dots, N$. Furthermore, we discretize (\ref{eq:NuRe:3}) using a standard central finite differences scheme to obtain
\begin{equation} \label{eq:NuRe:5}
	\frac{d}{dt} \mathbf z = \mathbb{J}_{2N} L\mathbf z,
\end{equation}
where $\mathbf z=(q,\dots,q_N,p_q,\dots,p_n)^T$ and
\begin{equation} \label{eq:NuRe:6}
L = 
\begin{pmatrix}
	I_n & 0_N \\
	0_N & \kappa(\omega)D_{xx}
\end{pmatrix},\quad 
\end{equation}
with $D_{xx}$ the central finite differences matrix operator. The discrete Hamiltonian can finally be written as
\begin{equation} \label{eq:NuRe:7}
	H_{\Delta x}(\mathbf z) = \frac{\Delta x}2 \sum_{i=1}^{N} \left( p_i^2 + \kappa(\omega) \frac{(q_{i+1} - q_i)^2}{2\Delta x ^ 2} + \kappa(\omega) \frac{(q_{i} - q_{i-1})^2}{2\Delta x ^ 2} \right).
\end{equation}
The initial condition is given by
\begin{equation} \label{eq:NuRe:8}
	q_i(0) = h( 10\times|x_i - \frac{1}{2}| ), \quad p_i = 0, \quad i=1,\dots,N
\end{equation}
where $h(s)$ is the cubic spline function
\begin{equation} \label{eq:NuRe:9}
h(s) = 
\left\{
\begin{aligned}
& 1 - \frac{3}{2}s^2 + \frac{3}{4}s^3, \quad & 0\leq s \leq 1, \\
& \frac{1}{4}(2-s)^3, & 1< s \leq 2, \\
& 0, & s > 2.
\end{aligned}
\right.
\end{equation}
This will result in waves propagating in both directions on the torus.

For numerical time integration we {\edit use} the Str\"omer-Verlet (\ref{eq:Hasy:13}) scheme, {\edit which is explicit since} the Hamiltonian is separable for the linear wave-equation. The full model uses the following parameter set \vspace{0.5cm}
\begin{center}
\begin{tabular}{|l|l|}
\hline
Domain length & $L = 1$ \\
No. grid points & $N = 500$ \\
Space discretization size & $\Delta x = 0.002$ \\
Time discretization size & $\Delta t = 0.01$ \\
Wave speed & $c^2 = 0.1$ \\
\hline
\end{tabular}
\end{center}
\vspace{0.5cm}
We compare the reduced system obtained by the greedy algorithm with the methods based on SVD. To generate snapshots, we discretize the parameter space $[0,1]^4$ into in total of $5^4$ equidistant grid points. For the SVD based methods and POD, snapshots are gathered in the snapshot matrices $S$, $S_{\text{combined}}$ and $S_{\text{complex}}$, respectively, and the SVD is performed to construct the reduced basis. The greedy method is applied following Algorithm \ref{alg:SyMo:3}; as input, the tolerance for the error in the Hamiltonian is set to $\delta = 5 \times 10^{-3}$. All reduced systems are taken to have an identical size ($k=80$ for POD and $k=40$ for the symplectic methods). We use the {\edit Str\"omer-Verlet} scheme for symplectic methods and a second order Runge-Kutta method for the POD. {\edit The choice of different time integration routines is due to the fact that the POD destroys the canonical form of the original equations and a symplectic integrator cannot be applied. One can alternatively use separate reduced subspaces for the potential and the momentum spaces, which however is not a standard model reduction approach and requires further analysis.} Finally we use transformation (\ref{eq:SyMo:1}) to transfer the solution of the reduced systems into the high-dimensional space for illustration purposes.

We reduced the cost by 50\% in the offline stage when using the greedy method as compared to SVD-based methods (cotangent lift and complex SVD method). This happens because the SVD-based methods require time integration of the full system for all discrete parameter points, while the greedy method picks a number of parameters from the parameter space.

%
%

\begin{figure}

\begin{minipage}{.5\linewidth}
\centering
\subfloat[$t=0$]{\label{fig:NuRe:1a}\includegraphics[width=1\textwidth]{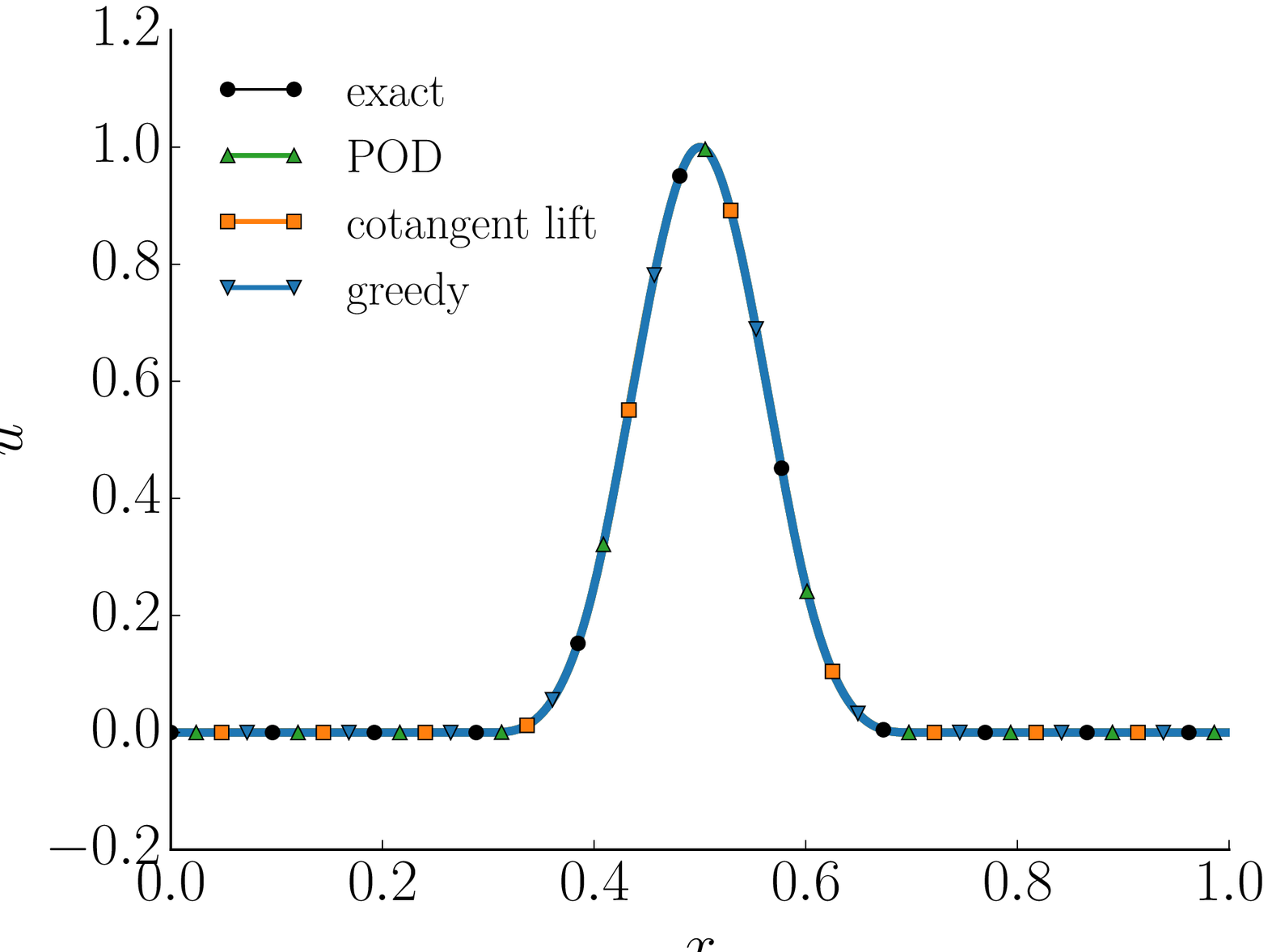}}
\end{minipage}%
\begin{minipage}{.5\linewidth}
\centering
\subfloat[$t=1$]{\label{fig:NuRe:1b}\includegraphics[width=\textwidth]{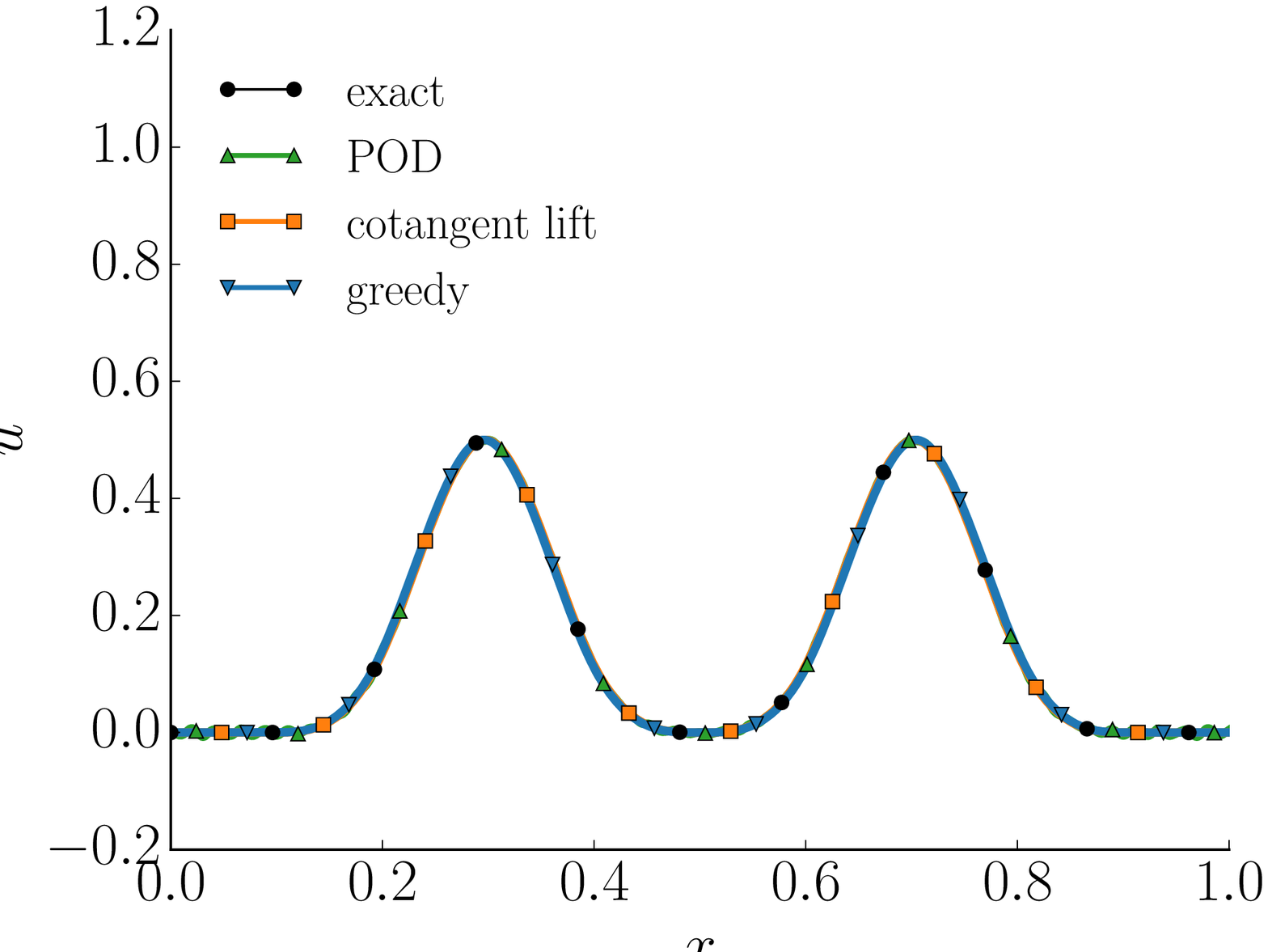}}
\end{minipage}\par\medskip
\centering
\subfloat[$t=2$]{\label{fig:NuRe:1c}\includegraphics[width=0.5\textwidth]{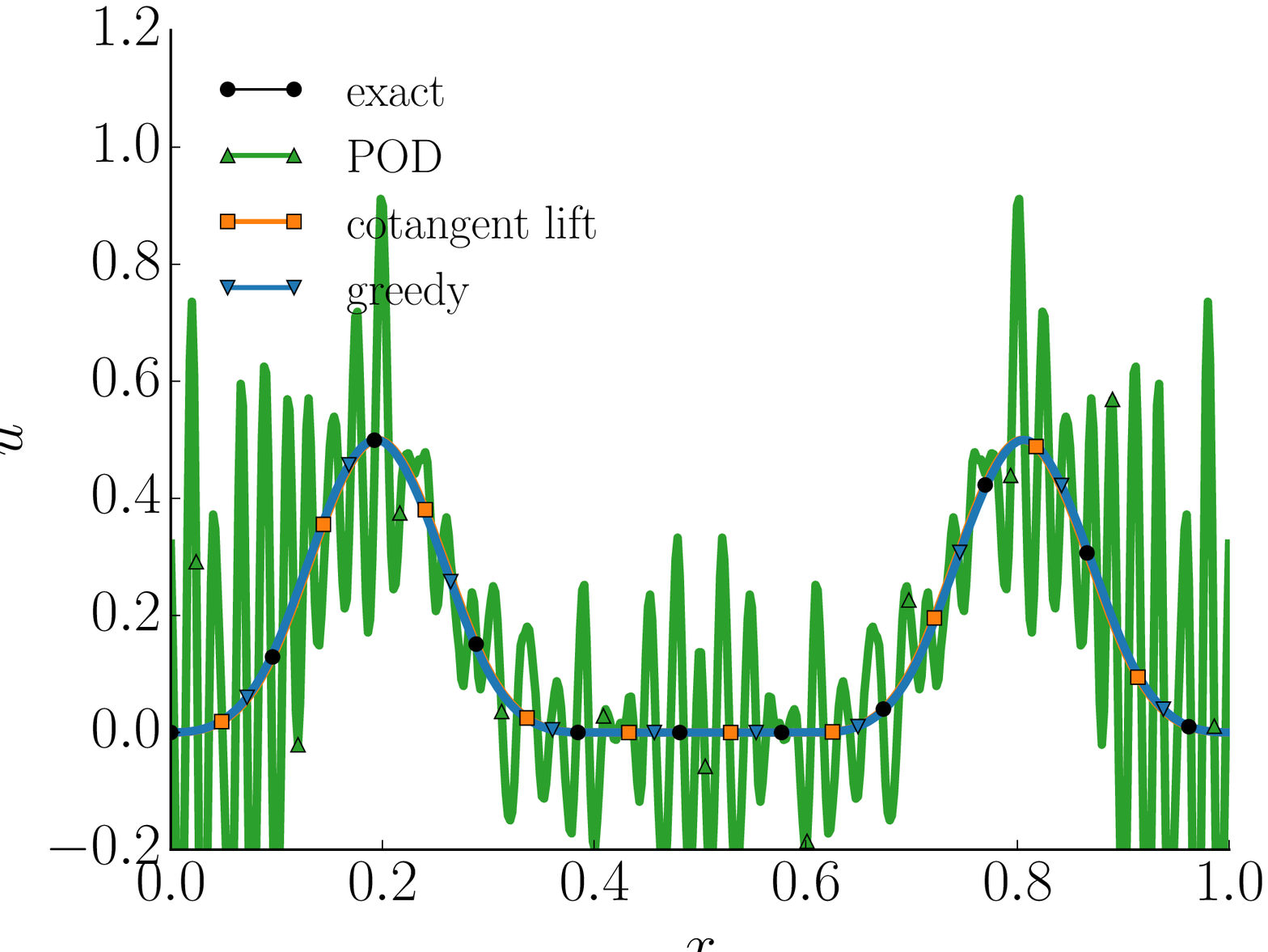}}
\caption{The solution $q$ at $t=0$, $t=1$ and $t=2$ of the linear wave equation for parameter value $c= 0.1019$ different from training parameters. Here, the solution of the full system together with the solution of the POD, cotangent lift, complex SVD and the greedy reduced system is shown.}
\label{fig:NuRe:1}
\end{figure}

\begin{figure}

\begin{minipage}{.5\linewidth}
\centering
\subfloat[]{\label{fig:NuRe:2c}\includegraphics[width=\textwidth]{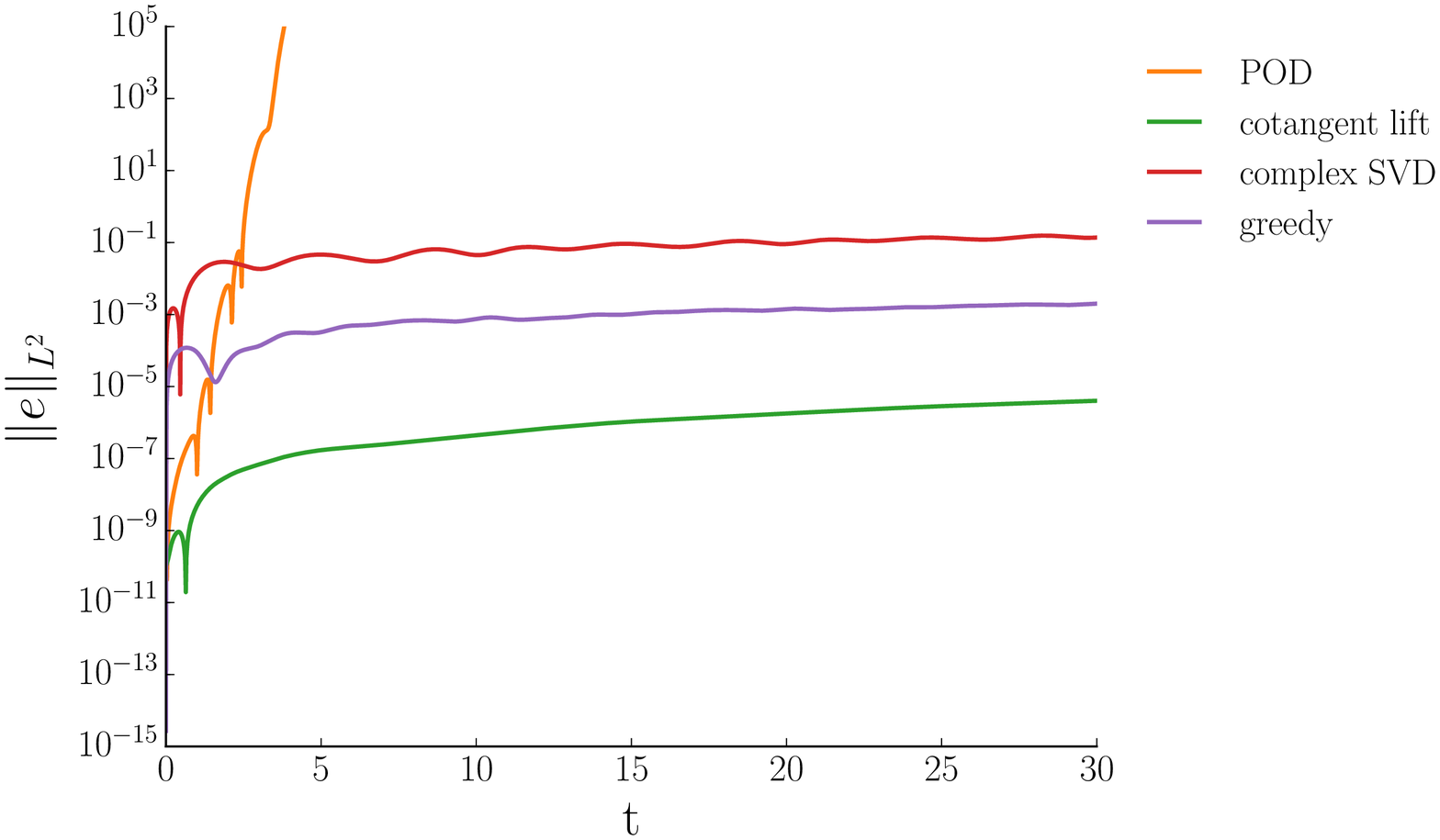}}
\end{minipage}%
\begin{minipage}{.5\linewidth}
\centering
\subfloat[]{\label{fig:NuRe:2b}\includegraphics[width=\textwidth]{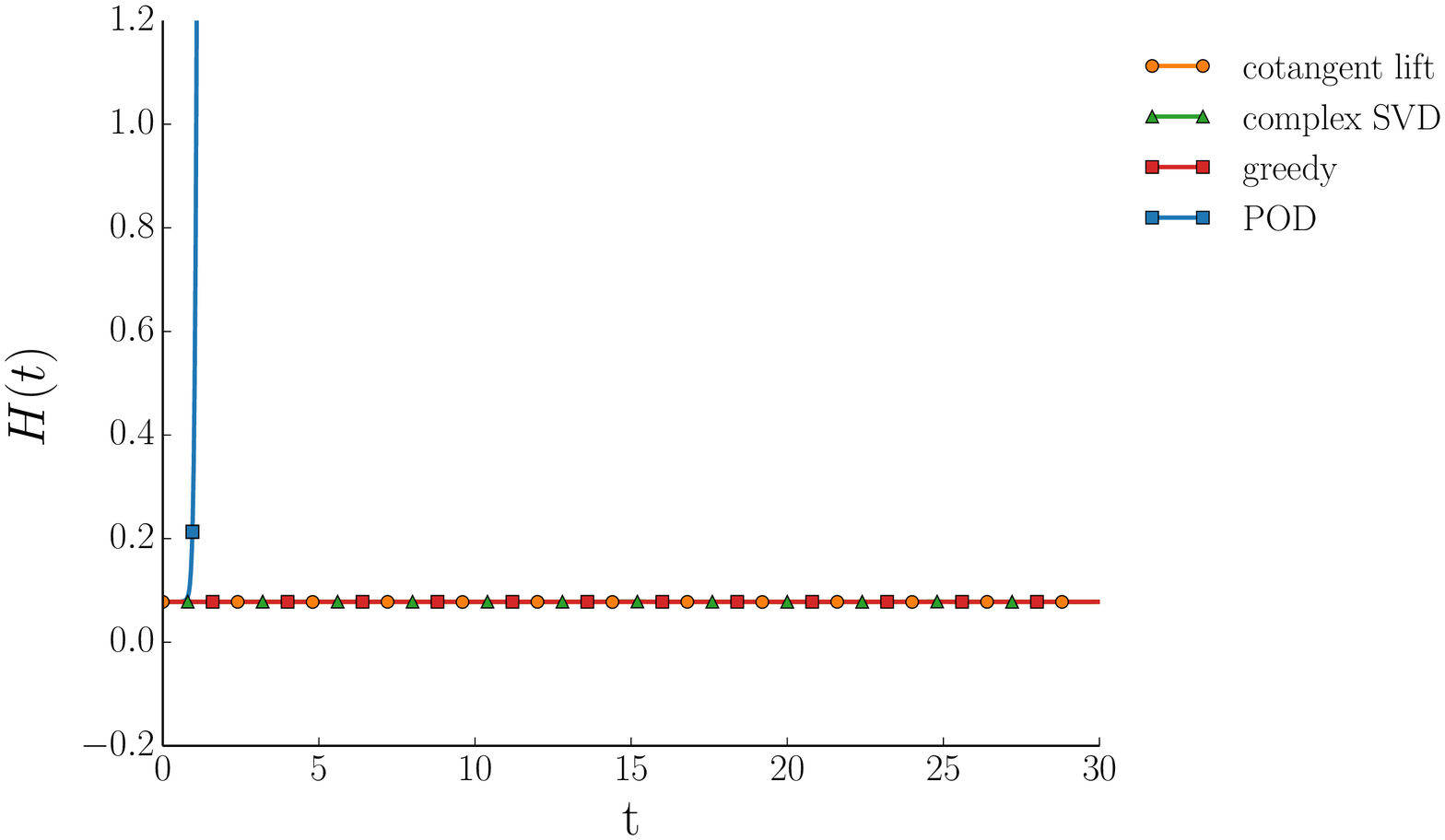}}
\end{minipage}\par\medskip
\centering

\caption{{\edit (a) The $L^2$-error between the solution of the full system and the reduced system for different model reduction methods for $t \in [0,30]$.} (b) Plot of the Hamiltonian function for $t \in [0,30]$. }
\label{fig:NuRe:2}
\end{figure}

Figure \ref{fig:NuRe:1a} shows the solution of the linear wave equation for parameter values {\edit $(\omega_1,\omega_2,\omega_3,\omega_4) = (0.8456,0.1320,0.9328,0.5809)$ or $\kappa(\omega) = 0.1019$}, chosen to be different from training parameters, at $t=0$, $t = 1$ and $t=2$. While we see instability and divergence from the exact solution for the POD reduced system, the symplectic methods provide a good approximation of the full model. 

{\edit The decay of the singular values for the POD are shown in Figure \ref{fig:NuRe:5a}. The decay} of the singular values suggests that a low dimensional solution manifold indeed exists. However, since the linear subspace, constructed by the POD, is not symplectic, we observe blow up of the Hamiltonian function in Figure \ref{fig:NuRe:2b} and the instability of the solution in Figure \ref{fig:NuRe:1}. The symplectic methods (using a reduced basis of the same size as POD) preserve the Hamiltonian function as shown in Figure \ref{fig:NuRe:2b}.

Figure \ref{fig:NuRe:2c} shows the $L^2$-error between the solution of the full model and the reduced systems constructed by different methods. We note that the error for the POD reduced system rapidly increases, confirming that {\edit the} projection based reduced system does not yield a stable solution. Furthermore, the symplectic methods provide a better approximation since the geometric structure of the original system is preserved. Although the greedy method is almost {\edit twice faster than the SVD-based methods} in the offline stage, its accuracy is comparable. {\edit The cotangent lift method provides a more accurate solution, on the other hand the cotangent lift basis (\ref{eq:SyMo:9}) takes a less general form and usually computationally more demanding than the greedy method.}

{\edit For complex systems were the solution of the full system is expensive and for high dimensional parameter domains, POD-based methods become impractical \cite{Anonymous:2016wl,Quarteroni:2016wi}. However, the greedy method requires substantially fewer (proportional to the size of the reduced basis) evaluation of the time integration of the original system.}






\subsection{Nonlinear Schr\"odinger equation} \label{chap:NuRe:1.2} Let us consider the one-dimensional parametric Schr\"odinger equation
\begin{equation} \label{eq:NuRe:10}
\left\{
\begin{aligned}
	& i u_t(t,x,\epsilon) = - u_{xx}(t,x,\epsilon) - \epsilon |u(t,x,\epsilon)|^2 u(t,x,\epsilon),\\
	& u(0,x) = u_0(x),
\end{aligned}
\right.
\end{equation}
where $u$ is a complex valued wave function, $i$ is the imaginary unit, $|\cdot|$ is the modulus operator and $\epsilon$ is a parameter that belongs to the interval $\Gamma = [0.9,1.1]$. We consider periodic boundary conditions, i.e., $x$ belongs to a one-dimensional torus of length $L$. We consider the initial condition
\begin{equation} \label{eq:NuRe:11}
	u_0(x) = \frac{\sqrt 2}{\cosh(x - x_0)} \exp(i\frac{c(x-x_0)}{2}),
\end{equation}
for a positive constant $c$. In quantum mechanics, the quantity $|u(t,x)|^2$ represents the probability of finding the system in state $x$ at time $t$. For the choice of $\epsilon = 1$, $|u(x,t)|$ becomes a solitary wave, and the initial condition will be transported in the positive $x$ direction with a constant speed. For other choices of $\epsilon$, the solution comprises an ensemble of solitary waves, moving in either direction \cite{Faou:2012vh}. 

By introducing the real and imaginary variables $u = p + iq$, we can rewrite (\ref{eq:NuRe:10}) in canonical form as
\begin{equation} \label{eq:NuRe:12}
\left\{
\begin{aligned}
 q_t &= p_{xx} + \epsilon (q^2+p^2)p, \\
 p_t &= -q_{xx} - \epsilon (q^2 + p^2)q,
\end{aligned}
\right.
\end{equation}
with the Hamiltonian function
\begin{equation} \label{eq:NuRe:13}
	H(q,p) = \int_{0}^{L} (q_x^2 + p_x^2) + \frac \epsilon 2 (q^2 + p^2)^2\ dx.
\end{equation}
We discretize the torus into $N$ equidistant points and take $\Delta x = L/N$, $x_i = i\Delta x$, $q_i=q(t,x_i,\epsilon)$ and $p_i = p(t,x_i,\omega)$ for $i = 1 ,\dots,N$. A central finite differences scheme is used to discretize (\ref{eq:NuRe:12}) as
\begin{equation}  \label{eq:NuRe:14}
	\frac{d}{dt} \mathbf z = \mathbb J_{2N} L\mathbf z + \mathbb J_{2N} \mathbf g(\mathbf z).
\end{equation}
Here $\mathbf z = (q_1,\dots,q_N,p_1,\dots,p_n)^T$ and
\begin{equation}  \label{eq:NuRe:15}
	L = 
	\begin{pmatrix}
		D_{xx} & 0_N \\
		0_N & D_{xx}
	\end{pmatrix}.
\end{equation}
Here $\mathbf g$ is a vector valued nonlinear function defined as
\begin{equation}  \label{eq:NuRe:16}
	\mathbf g(\mathbf z) =
	\begin{pmatrix}
	(q_1^2 + p_1^2)q_1 \\
	\vdots \\
	(q_N^2 + p_N^2)q_N \\
	(q_1^2 + p_1^2)p_1 \\
	\vdots \\
	(q_N^2 + p_N^2)p_N
	\end{pmatrix}.
\end{equation}
We discretize the Hamiltonian to obtain
\begin{equation}  \label{eq:NuRe:17}
	H_{\Delta x}(\mathbf z) = {\Delta x}\sum_{i=1}^{N} \left( \frac{q_i q_{i-1} - q_i^2}{\Delta x ^2} + \frac{p_i p_{i-1} - p_i^2}{\Delta x ^2} + \frac \epsilon 4 (p_i^2 + q_i^2)^2  \right),
\end{equation}
and use a Str\"omer-Verlet (\ref{eq:Hasy:13}) scheme for time integration. Since the Hamiltonian function (\ref{eq:NuRe:17}) is non-separable, this scheme becomes implicit so in each time iteration, a system of nonlinear equations is solved using Newton's iteration. We summarize the physical and numerical parameters for the full model in the following table

\vspace{0.5cm}
\begin{center}
\begin{tabular}{|l|l|}
\hline
Domain length & $L = 2\pi /l$ \\
Domain scaling factor & $l = 0.11$ \\
wave speed & $c =1$\\
No. grid points & $N = 256$ \\
Space discretization size & $\Delta x = 0.2231$ \\
Time discretization size & $\Delta t = 0.01$ \\
\hline
\end{tabular}
\end{center}
\vspace{0.5cm}
Regarding computation of the nonlinear terms of reduced systems, {\edit we} compare the DEIM with the symplectic DEIM. For generation of the DEIM reduced basis we apply Algorithm \ref{alg:MoOr:1} to the set of nonlinear snapshots. Algorithm \ref{alg:SyMo:4} is used to construct a reduced basis appropriate for the symplectic DEIM. As input, we provide the symplectic basis generated by Algorithm \ref{alg:SyMo:3} with the set of nonlinear snapshots and a tolerance for the error $\delta = 10^{-4}$.

We compare the reduced system obtained using the greedy algorithm with the cotangent lift, the complex SVD, DEIM, the symplectic DEIM and also the POD. For the SVD-based methods, we discretize the parameter space $[0.9,1.1]$ into $M=500$ equidistant grid points {\edit across} the discrete parameter space $\Gamma_M = \{\epsilon_1,\dots,\epsilon_M \}$, and gather trajectory snapshots for each $\epsilon_i$ for $i = 1,\dots,M$ in the snapshots matrix $S$. All reduced systems are taken to have identical sizes ($k=90$ for the symplectic methods and $k=180$ for the POD method). Following Algorithm \ref{alg:SyMo:3} we construct the reduced system using the same discrete parameter space $\Gamma_M$. The tolerance for the error in the Hamiltonian is set to $\delta = 10^{-3}$. Moreover, for DEIM and symplectic DEIM, we construct bases of size $k'=80$. Note that the reduced system, generated in the symplectic DEIM, will be of size $k+k'=170$.

The cost of the offline stage is reduced to 20\% when using the greedy method for constructing a symplectic basis of size $k=90$, as compared to the SVD-based methods. The online stage, i.e., time integration for a new parameter in $\Gamma$, is generally more than 3 times faster than {\edit for} the original system. We point out that the efficiency of reduced systems are implementation and platform dependent {\edit and we expect further reduction as the size of the problem increases.}

\begin{figure}

\begin{minipage}{.5\linewidth}
\centering
\subfloat[$t=0$]{\label{fig:NuRe:3a}\includegraphics[width=1\textwidth]{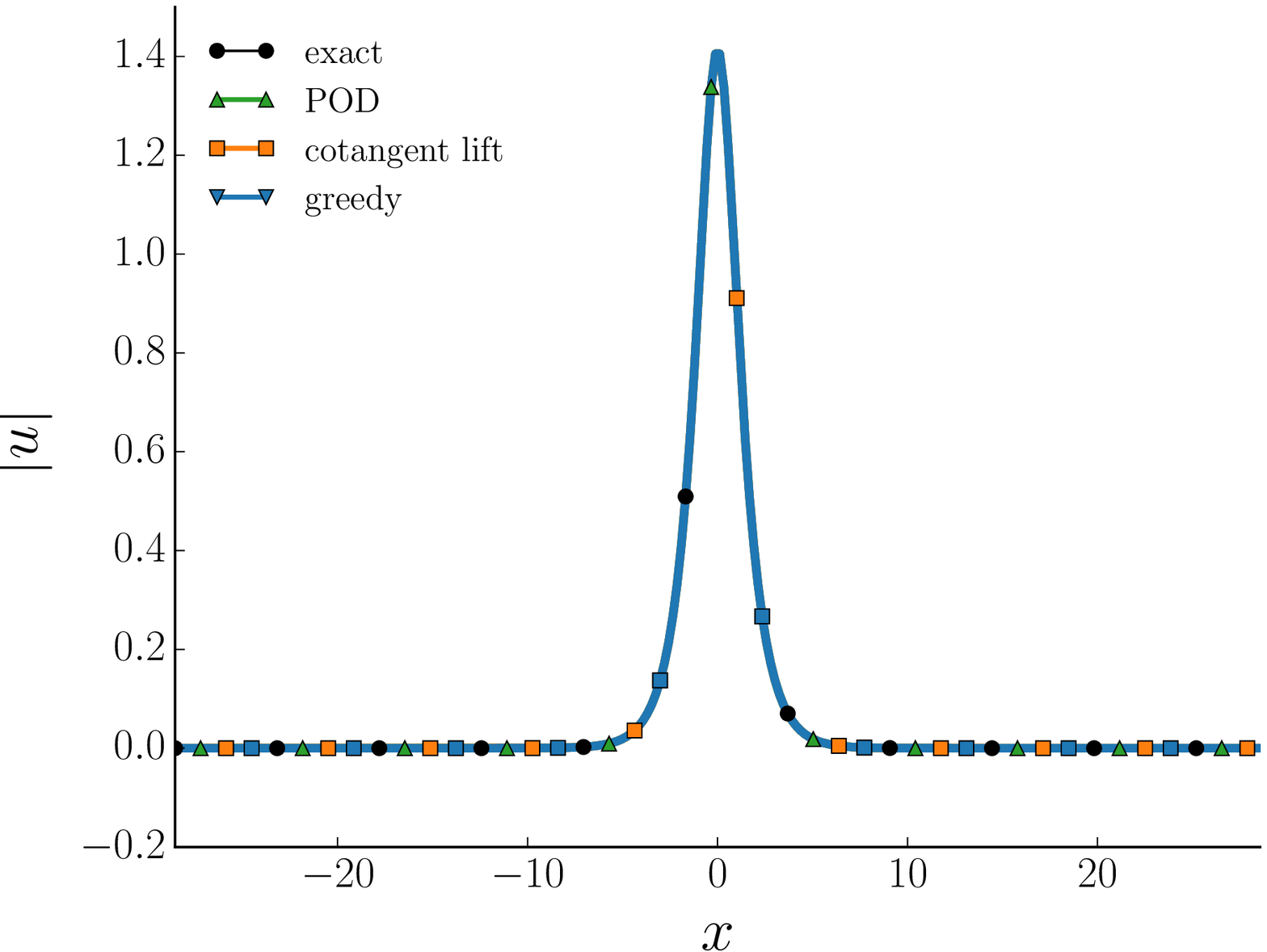}}
\end{minipage}%
\begin{minipage}{.5\linewidth}
\centering
\subfloat[$t=10$]{\label{fig:NuRe:3b}\includegraphics[width=\textwidth]{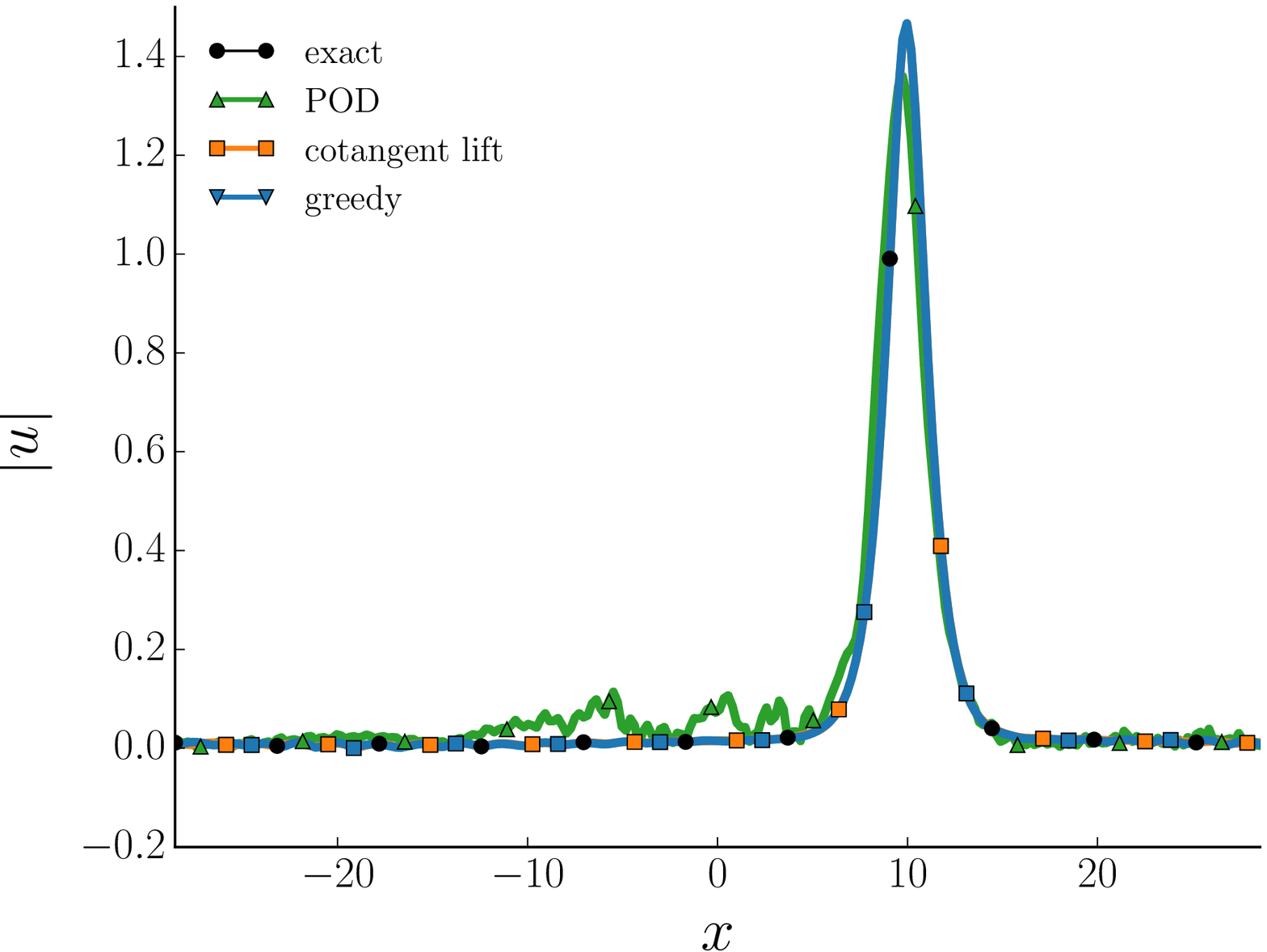}}
\end{minipage}\par\medskip
\centering
\subfloat[$t=20$]{\label{fig:NuRe:3c}\includegraphics[width=0.5\textwidth]{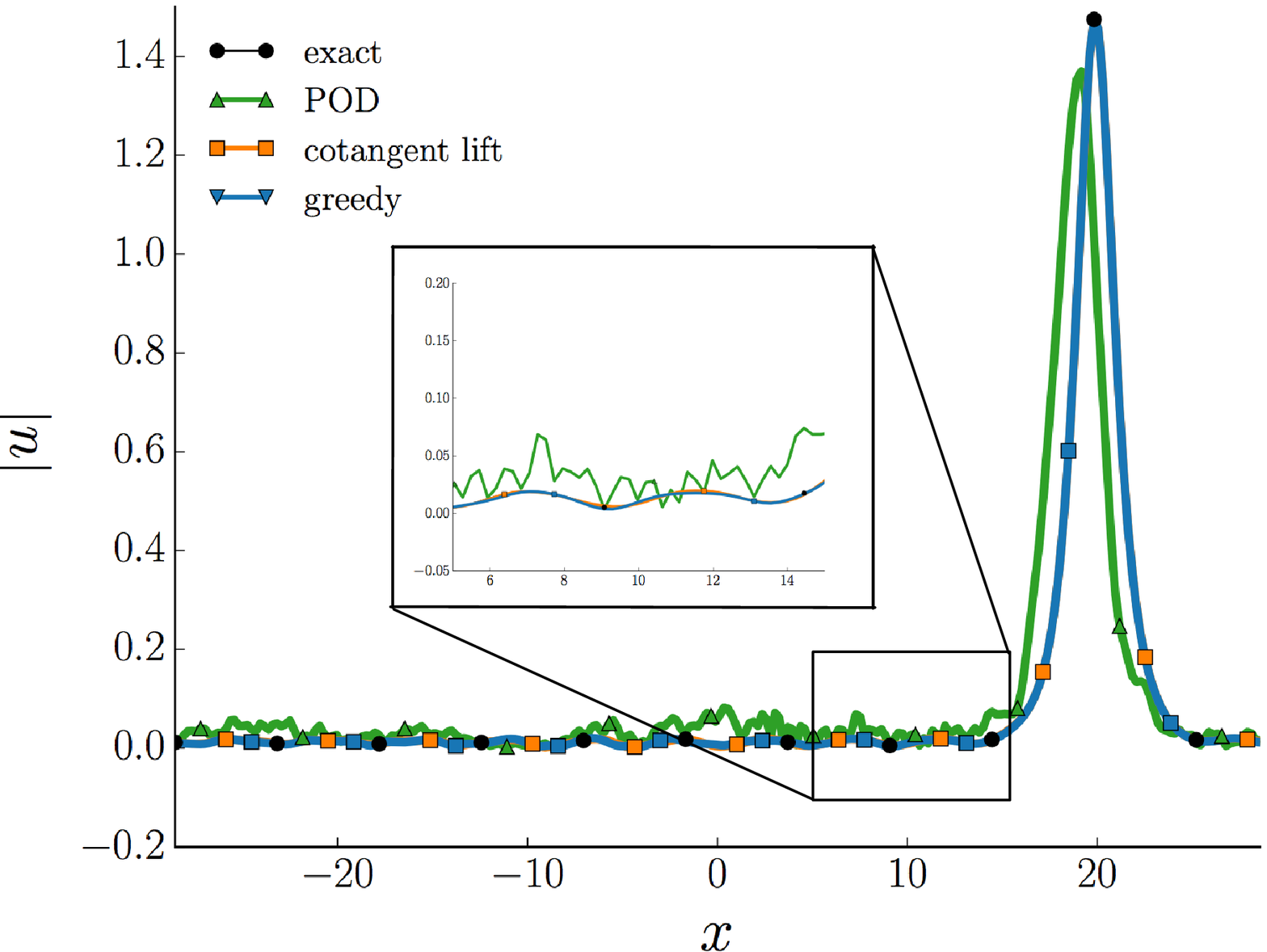}}
\caption{The solution $|u(t,x)| = \sqrt{q^2 + p^2}$ at $t=0$, $t=10$ and $t=20$ of the Nonlinear Schr\"odinger equation for parameter value $\epsilon = 1.0932$. Here the solution of {\edit the} full system, together with the solution of the POD, cotangent lift, complex SVD and the greedy reduced system, is shown.}
\label{fig:NuRe:3}
\end{figure}

\begin{figure}

\begin{minipage}{.5\linewidth}
\centering
\subfloat[]{\label{fig:NuRe:4c}\includegraphics[width=\textwidth]{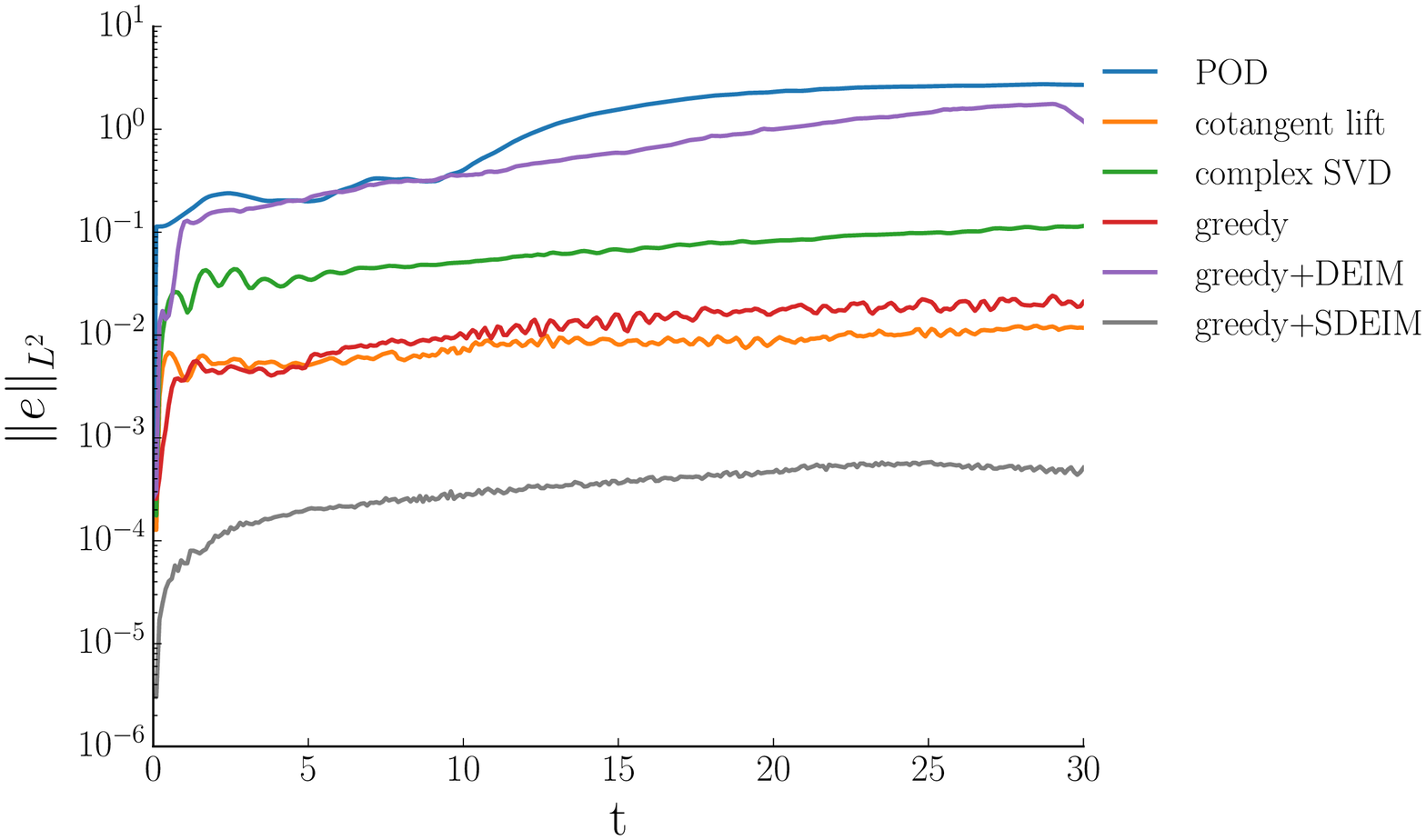}}
\end{minipage}%
\begin{minipage}{.5\linewidth}
\centering
\subfloat[]{\label{fig:NuRe:4b}\includegraphics[width=\textwidth]{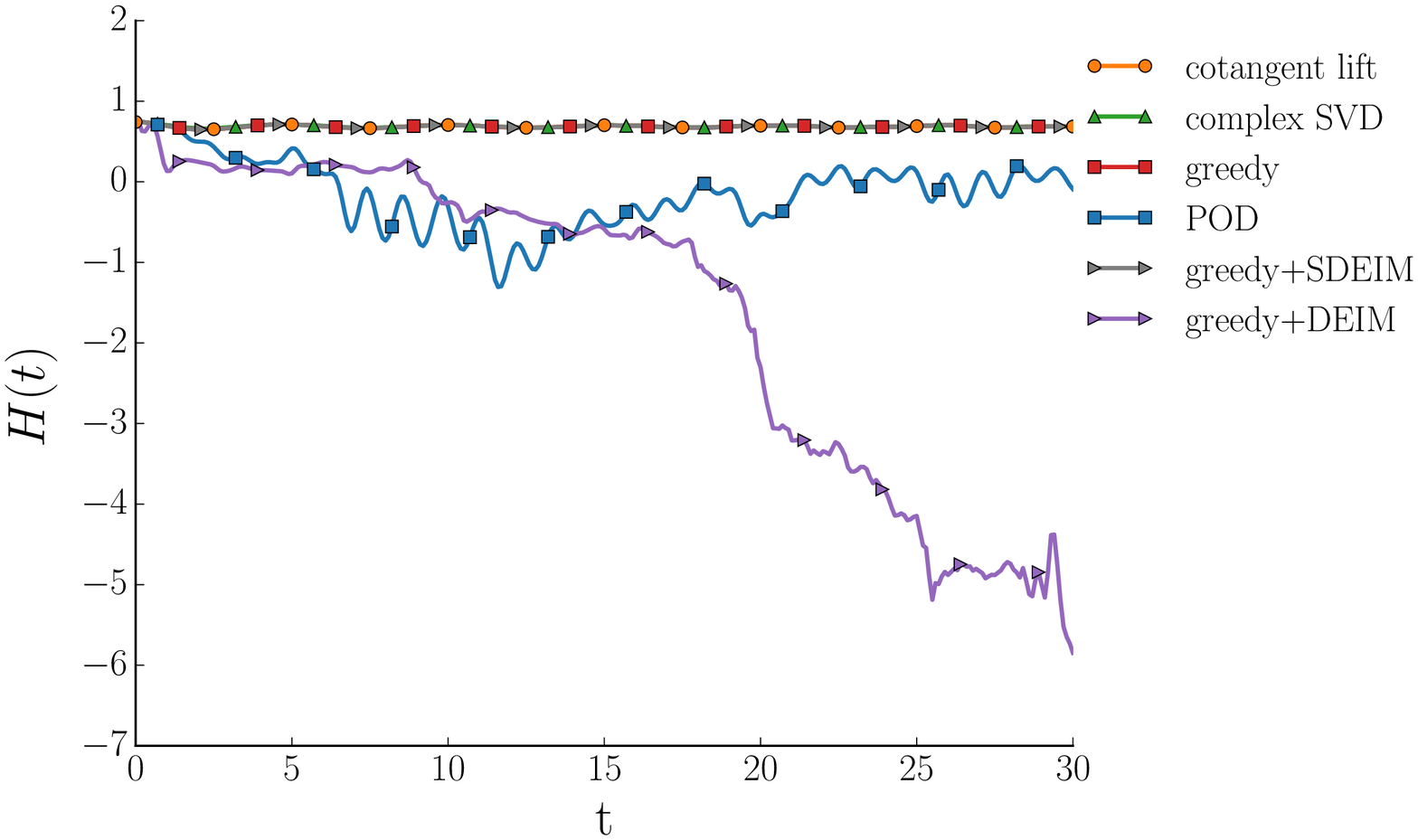}}
\end{minipage}\par\medskip
\centering

\caption{{\edit (a) Plot of the Hamiltonian function for $t \in [0,30]$.} (b) The $L^2$ error between the solution of the full system and the reduced system for different model reduction methods for $t \in [0,30]$. }
\label{fig:NuRe:4}
\end{figure}

Figure \ref{fig:NuRe:3} shows the solution of the Schr\"odinger equation for parameter value $\epsilon = 1.0932$ at $t=0$, $t=10$ and $t=20$. We first compare the reduced system obtained by the greedy algorithm with the POD, the cotangent lift, and the complex SVD method. The size of the reduced systems are taken identical for all methods ($k=180$ for POD and $k=90$ for the rest). Although the decay of the singular values in Figure \ref{fig:NuRe:5b} suggests that the accuracy of the POD reduced system should be comparable to that of the other methods, we observe instabilities in the solution at $t=10$. The greedy, the cotangent lift and the complex SVD method, on the other hand, generate a stable reduced system that accurately approximates the solution of the full model.

{\edit In Figure \ref{fig:NuRe:4b} we observe that the symplectic methods preserve the Hamiltonian function, unlike the POD and the DEIM methods. We emphasise that using the reduced basis, obtained by the greedy, together with the DEIM (purple line) does not preserve the symplectic structure as suggested in this figure.}

Figure \ref{fig:NuRe:4c} illustrates the $L^2$-error between the solution of the full model with the reduced systems, generated by different methods. We first observe that symplectic methods yield a lower computational error {\edit when} compared to non-symplectic methods. Secondly, we observe that although the reduced systems from the cotangent lift and the complex SVD are of the same size, their accuracy is different by an order of magnitude. We notice that the greedy algorithm is slightly less accurate than the cotangent lift method while its offline computational cost is reduced to 20\% {\edit when compared to the cotangent lift}. Lastly we notice that the combination of the greedy reduced basis and DEIM yields large errors in the solution while the solution using the symplectic DEIM is very accurate. We note that the symplectic DEIM is even more accurate than the greedy itself since it has been enriched by the nonlinear snapshots. 

\subsection{Numerical Convergence} 
\begin{figure}

\begin{minipage}{.5\linewidth}
\centering
\subfloat[]{\label{fig:NuRe:5a}\includegraphics[width=\textwidth]{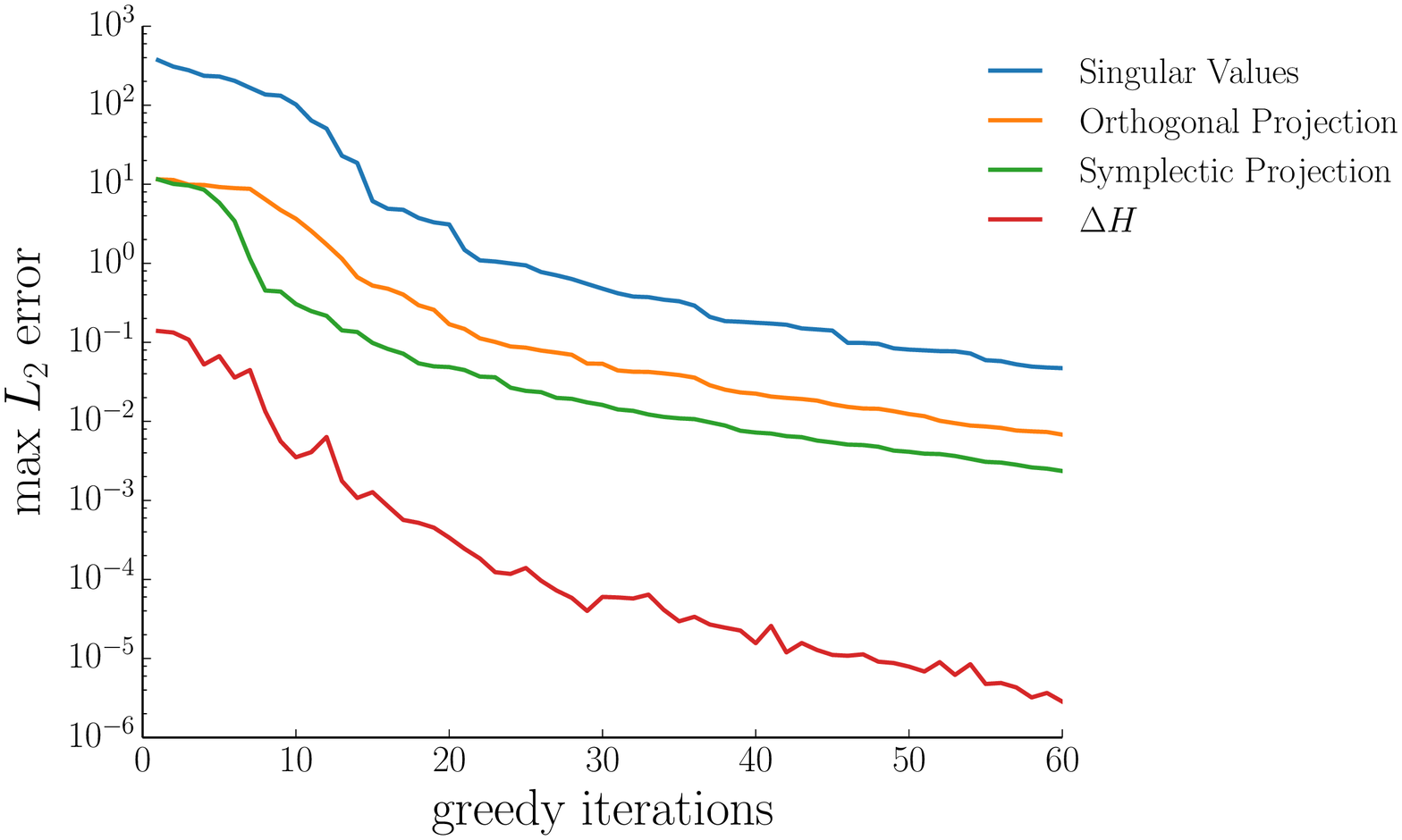}}
\end{minipage}%
\begin{minipage}{.5\linewidth}
\centering
\subfloat[]{\label{fig:NuRe:5b}\includegraphics[width=\textwidth]{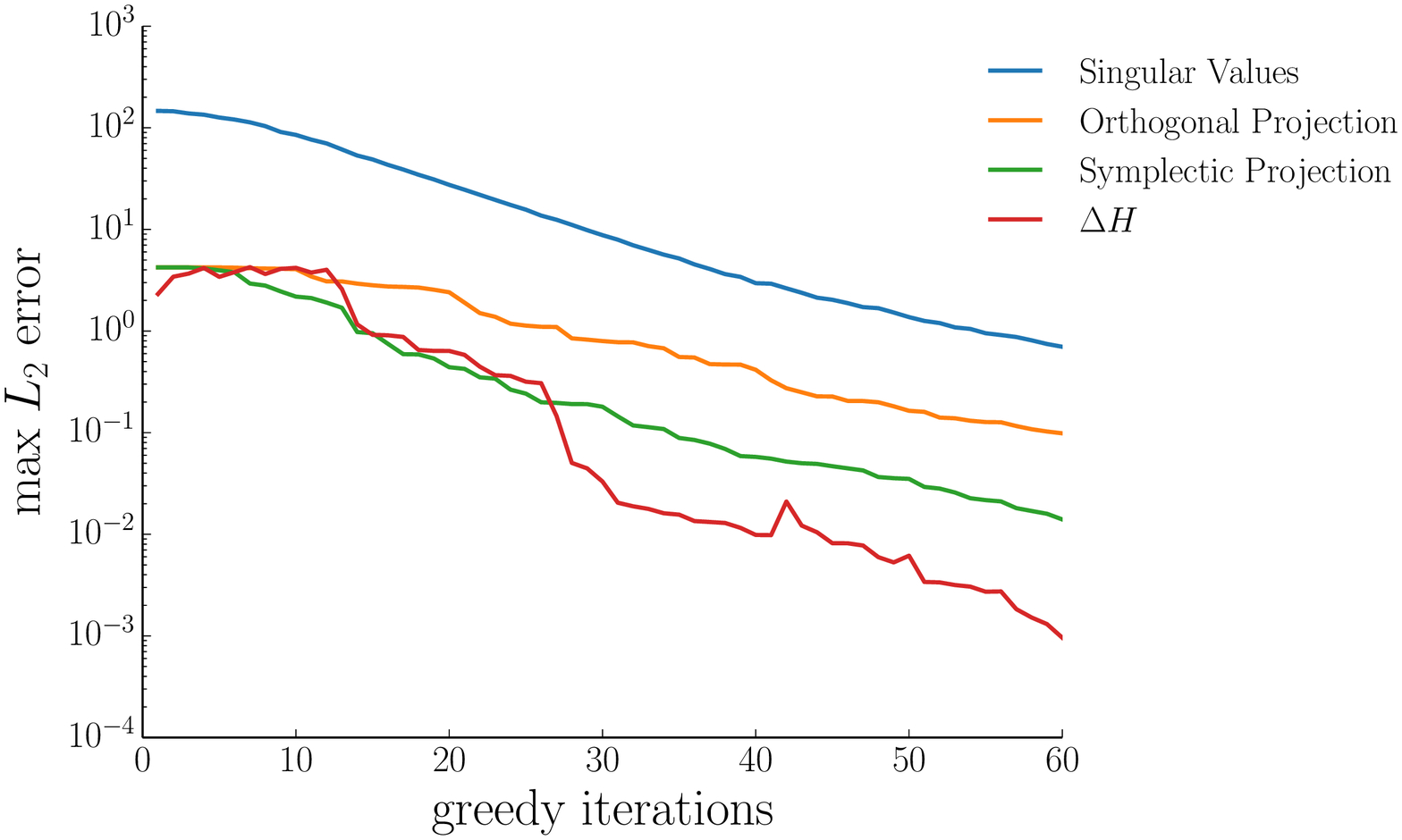}}
\end{minipage}\par\medskip
\centering

\caption{ (a) Convergence of the greedy method for the wave equation. (b) Convergence of the greedy method for the nonlinear Schr\"odinger equation equation. }
\label{fig:NuRe:5}
\end{figure}

{\edit In this section we discuss the numerical convergence of the symplectic greedy method introduced in Section \ref{chap:SyMo:1}. The exponential convergence properties of the conventional greedy \cite{Quarteroni:2016wi} is presented in \cite{Buffa:2012iz,Binev:2011fj}. Theorem \ref{theorem:SyMo:2} suggests that the symplectic greedy method has similar properties. To illustrate this we compare the convergence of the conventional greedy with the convergence of the symplectic greedy method through the numerical simulations in Sections \ref{chap:NuRe:1.1} and \ref{chap:NuRe:1.2}. 

The decay of the singular values of the snapshot matrix for the parametric wave equation and the nonlinear Schr\"odinger equation are given in Figure \ref{fig:NuRe:5}. The decay rate of the singular values is a strong indicator for the decay rate of the Kolmogorov $n$-width of the solution manifold. We expect that the conventional greedy method and the symplectic greedy method provide a similar rate in the decay of the error.
	
Figure \ref{fig:NuRe:5} shows the maximum $L^2$ error between the original system and the reduced system at each iteration of different greedy methods. In this figure we find the conventional greedy with orthogonal projection error as a basis selection criterion (orange), the symplectic greedy method with a symplectic projection error as a basis selection criterion (green), and the symplectic greedy method with energy loss $\Delta H$ as a basis selection criterion (red).

It is observed that the decay rate of the error for greedy with the orthogonal projection and the greedy with the symplectic projection is similar to the decay of the singular values. This matches our expectation from Theorem \ref{theorem:SyMo:2}. We also notice that the greedy method with the loss in Hamiltonian provides an excellent error indication as a basis selection criterion.
}

\section{Conclusion} \label{chap:Con:1}

In this paper, we present a greedy approach for the construction of a reduced system that preserves the geometric structure of Hamiltonian systems. An iteration of the greedy method comprises searching the parameter space using the error in the Hamiltonian, to find the best basis vectors that increase the overall accuracy of the reduced basis. We argue that for a compact subset with exponentially small Kolmogorov $n$-width we recover exponentially fast convergence of the greedy algorithm. For fast approximation of nonlinear terms, the basis obtained by the greedy was combined with a symplectic DEIM to construct {\edit a reduced system with a Hamiltonian that is arbitrary close to the Hamiltonian of the original system.}


The numerical results demonstrate that the greedy method can save substantial computational cost in the offline stage as compared to alternative SVD-based techniques. Also since the reduced system obtained by the greedy method is Hamiltonian, the greedy method yields a stable reduced system. Symplectic DEIM effectively reduces computational cost of approximating nonlinear terms while preserving stability and symplectic structure. Hence, the greedy method is an efficient model reduction technique that provides an accurate and stable reduced system for large-scale parametric Hamiltonian systems.


\section*{Acknowledgments}
{\edit We would like to thank the referees for providing us with very useful comments which served to improve the paper.}


\bibliographystyle{siamplain}
\bibliography{ref}
\end{document}